\documentclass[letterpaper, 10 pt, conference]{ieeeconf}  

\IEEEoverridecommandlockouts                              

\usepackage{graphics} 
\usepackage{epsfig} 
\usepackage{amsmath} 
\usepackage{amssymb}  
\newcommand\numberthis{\addtocounter{equation}{1}\tag{\theequation}}

\usepackage{xcolor} 
\usepackage{pgfplots} 
\newtheorem{theorem}{Theorem}
\newtheorem{remark}{Remark}
\newtheorem{lemma}[theorem]{Lemma}

\usepackage{tikz}
\usetikzlibrary{shapes.geometric,arrows}
\usepackage[]{todonotes}

\usepackage{fancyhdr,lipsum}
\fancypagestyle{firstpage}{
  \fancyhf{}
  \fancyhead[C]{This is an author's copy supplementing a paper submitted to ACC 2022 and will be copyrighted by IEEE if accepted.}
}

\title{\LARGE \bf Robust Performance Analysis of Source-Seeking Dynamics with Integral Quadratic Constraints}

\author{Adwait Datar$^{1}$ and Herbert Werner$^{1}$
\thanks{This work was funded by the German Research Foundation (DFG) within
	their priority programme SPP 1914 Cyber-Physical Networking.}
\thanks{$^{1}$Institute of Control systems, Hamburg University of Technology,
	Ei{\ss}endorfer Str. 40, 21073 Hamburg, Germany.
       {\tt\small \{adwait.datar, h.werner\}@tuhh.de}}%
}

\pgfplotsset{compat=1.17}
\begin{document}

\maketitle
\thispagestyle{firstpage}
\pagestyle{empty}

\begin{abstract}
We analyze the performance of source-seeking dynamics involving either a single vehicle or multiple flocking-vehicles embedded in an underlying strongly convex scalar field with gradient based forcing terms.
For multiple vehicles under flocking dynamics embedded in quadratic fields, we show that the dynamics of the center of mass are equivalent to the dynamics of a single agent.
We leverage the recently developed framework of $\alpha$-integral quadratic constraints (IQCs) to obtain convergence rate estimates.
We first present a derivation of \textit{hard} Zames-Falb (ZF) $\alpha$-IQCs involving general non-causal multipliers based on purely time-domain arguments and show that a parameterization of the ZF multiplier, suggested in the literature for the standard version of the ZF IQCs, can be adapted to the $\alpha$-IQCs setting to obtain quasi-convex programs for estimating convergence rates.
Owing to the time-domain arguments, we can seamlessly extend these results to linear parameter varying (LPV) vehicles possibly opening the doors to non-linear vehicle models with quasi-LPV representations. 
We illustrate the theoretical results on a linear time invariant (LTI) model of a quadrotor, a non-minimum phase LTI plant and two LPV examples which show a clear benefit of using general non-causal dynamic multipliers to drastically reduce conservatism.
\end{abstract}

\section{Introduction}
Source-seeking problems are motivated by practical problems such as finding the source of an oil spill \cite{senga2007development}.
The abstract problem involves one or more vehicles located at arbitrary locations in an underlying scalar field with the goal of moving towards the minimum (or maximum based on convention) of the field which is called the source.
These vehicles are typically assumed to be able to estimate the gradient of the field at their respective locations either via communication with others or by making several measurements around their positions (see \cite{khong2014multi}. \cite{ogren2004cooperative}, \cite{Datar.51220205152020} for details).
Although the theoretical results are under the assumption of strongly convex fields and availability of gradients, experimental results \cite{Datar.51220205152020} on quasi-convex and non-convex fields with noisy measurements show that with a pre-filtering step, which involves fitting a strongly convex field with data collected in the local neighborhood, these protocols perform reasonably well for some scenarios.

Integral quadratic constraints (IQCs) \cite{megretski1997system} is a framework for analyzing the stability of an interconnection of two systems; a linear and stable "nominal" system and a possibly non-linear "uncertain" system.
The key idea is to capture the input-output properties of the uncertain system in the form of integral quadratic constraints and use these constraints to derive sufficient conditions on the nominal system that certify stability of the interconnection. 
The IQCs are further classified into \textit{soft}-IQCs and stricter \textit{hard}-IQCs depending on the limits of the integrals involved. 
See \cite{J.Veenman.2013},\cite{UlfJonsson.} for a tutorial.

An alternative to the theory of IQCs is the framework of dissipativity which allows one to analyze the stability of the interconnection by searching for non-increasing energy-like functions called storage functions. 
The results obtained from \textit{hard}-IQCs can be interpreted with dissipativity based arguments and vice-versa\footnote{The link between the \textit{soft}-IQCs and dissipativity has only recently been clearly understood \cite{Scherer.16052021}}.
This equivalence connects our earlier work \cite{Datar.51220205152020}, \cite{Attallah.2020} on analyzing stability of source-seeking dynamics with dissipativity based arguments to this paper where we use IQCs.

\subsection*{Related Work}\label{sec:related_work}
Stability analysis for source-seeking dynamics within the flocking framework is studied in \cite{Datar.51220205152020}, \cite{Attallah.2020} albeit without any performance guarantees.
Our main motivation in this paper is to extend these results to performance analysis. 
Moreover, the stability analysis in \cite{Datar.51220205152020}, \cite{Attallah.2020} involves dissipativity based arguments by constructing valid (non-increasing) storage functions from physical energy-like functions.
These physically motivated storage functions have a diagonal structure (such as the one used in \cite{OlfatiSaber.2006}, \cite{Datar.51220205152020}, \cite{Attallah.2020}).
The idea is to automate the search for storage functions (not necessarily diagonal) and include a less conservative stability analysis certificate along with performance guarantees.
As pointed out earlier, these dissipativity based results can be interpreted as the ones obtained from \textit{hard}-IQCs and the latter is the choice for the theoretical development in this paper.
The starting points for the theoretical work in this paper are \cite{Lessard.2016} and \cite{Hu.2016} where the exponential versions of the IQCs are introduced for systems in discrete-time ($\rho-$IQCs) and in continuous-time ($\alpha-$IQCs), respectively.
\cite{Hu.2016} introduces the \textit{soft} Zames-Falb (ZF) $\alpha$-IQCs corresponding to causal multipliers.
While \cite{Zhang.25022019} extends \cite{Lessard.2016} to less conservative non-causal ZF multipliers in the discrete-time setting, \cite{Freeman.62720186292018} presents the extension of \cite{Hu.2016} to non-causal multipliers in the continuous-time setting.
The theory developed in \cite{Freeman.62720186292018} is in a very general setting of Bochner spaces and covers Lemma \ref{theom:lemma_pq_ZF} in this paper.
Moreover, proof of Lemma 3 from \cite{Freeman.62720186292018} (which corresponds to Lemma \ref{theom:lemma_pq_ZF} here) is not available and we therefore present a self-contained proof for our setting in the appendix.
\cite{Freeman.62720186292018} focuses on the derivation of the IQCs and does not consider parameterizations of the multiplier to arrive at a quasi-convex optimization problem for performance analysis.
The present paper extends these results by considering a parameterization proposed in \cite{J.Veenman.2013} adapted to the $\alpha$-IQCs setting.
\cite{Zhang.25022019} compares different multiplier factorizations which include the discrete-time analogue of the parameterization in \cite{J.Veenman.2013}.
Furthermore, we consider vehicle models which have integral action and present results in a manner that do not require a loop transformation.
We extend the simple time-domain based derivation of the \textit{hard} ZF IQCs for causal multipliers presented in a side-bar in \cite{Scherer.16052021} to the $\alpha$-IQCs case covering non-causal multipliers and illustrating the key role of convexity. 
Although \cite{Zhang.25022019} and \cite{Freeman.62720186292018} present examples showing the benefit of non-causal multipliers in the discrete-time case, we present an example of a continuous-time system with an integrator demonstrating the benefit of non-causal multipliers over causal ones.
A closely related work is \cite{Fazlyab.2018} where causal and static multipliers are used to obtain non-exponential convergence rates when the dynamics are not exponentially stable.
Finally, previous works (such as \cite{Lessard.2016}, \cite{Hu.2016} and \cite{Freeman.62720186292018}) on the exponential version of IQCs present results for linear time invariant (LTI) systems.
Since we consider the \textit{hard}-IQCs with purely time-domain arguments, we show that it is rather straightforward to extend these results to linear parameter varying (LPV) \cite{shamma1992linear} systems as is done in \cite{pfifer2015robustness} for the standard IQCs.
This opens the doors for considering non-linear vehicle models with quasi-LPV representations.
The framework developed in this paper can be applied to extremum-seeking control problems \cite{michalowsky2016extremum} or problems involving realtime-optimization \cite{nelson2018integral}. 
Although the setup considered in \cite{nelson2018integral} is the similar as in this paper (except for the LPV extension), the focus there is on designing an optimizer and obtaining conditions for optimality and stability with static $\alpha-$IQCs.
\subsection*{Outline:}
After discussing the notation in section \ref{sec:notation} and the problem setup in section \ref{sec:problem_setup}, the main theoretical results are presented in section \ref{sec:theory} followed by the extension to LPV systems in section \ref{sec:lpv_ext}. Numerical results are presented in section \ref{sec:numerical_results} leading to conclusions and open directions in section \ref{sec:conclusions}.
\section{Notation}\label{sec:notation}
Let $\mathbb{R}$ denote the set of real numbers.
Let $\mathbb{S}^{n}$ denote the set of symmetric matrices of size $n$.
For any $X \in \mathbb{S}^n$, $X>(\geq)0$ means that the matrix $X$ is symmetric positive definite (semi-definite), and $X<(\leq)0$ denote that $-X$ is symmetric positive definite (semi-definite).
For a non-singular symmetric matrix $X> 0$, let the positive number $\textnormal{cond}(X)$ be the condition number of $X$. 
For any $x\in \mathbb{R}^n$, let $\textnormal{diag}(x)$ denote the diagonal matrix formed by placing the entries of $x$ along the diagonal. 
For block matrices, we use $*$ to denote required entries to make the matrix symmetric. 
Let $\mathbf{0}$ and $\mathbf{1}$ denote the vectors or matrices of all zeroes and ones of appropriate sizes, respectively. 
Let $I_d$ be the identity matrix of dimension $d$ and we remove the subscript $d$ when the dimension is clear from context. 
Let $\otimes$ represent the Kronecker product.
Let $\mathcal{S}(m,L)$ denote the set of continuously differentiable functions $f:\mathbb{R}^d\rightarrow \mathbb{R}$ which are strongly convex with parameter $m$, and have Lipschitz gradients with parameter $L$ for any $0<m \leq L$, i.e., 
\begin{equation*}
    m||y_1-y_2||^2 \leq  (\nabla f(y_1)-\nabla f(y_2))^T(y_1-y_2) \leq L ||y_1-y_2||^2.
\end{equation*}
The set of vector valued functions which are square-integrable over $[0,T]$ for any finite $T<\infty$ is denoted by $\mathcal{L}_{2e}[0,\infty)$.
For any  signal $q \in \mathcal{L}_{2e}[0,\infty)$, and a fixed $T\geq 0$ let us define the extension 
\begin{equation}\label{eq:extension_signal}
    \begin{split}
        q_T(t)&=
        \begin{cases}
        q(t), & \text{if}\ t\in [0,T], \\
        0, & \text{if}\ t\in \mathbb{R} \backslash [0,T].
        \end{cases}
    \end{split}
\end{equation}
Let $\mathcal{L}_1(-\infty,\infty)$ denote the set of functions $h:\mathbb{R}\rightarrow \mathbb{R}$, such that $\int_{-\infty}^{\infty}|h(t)|dt < \infty$.
We use 
$\left[\begin{array}{c|c}
\mathcal{A}     &  \mathcal{B}\\
\hline
\mathcal{C}     &  \mathcal{D}
\end{array}\right]$
to represent an LTI system with state-space realization given by matrices $A,B,C$ and $D$.
The series interconnection of the LTI system 
$\left[\begin{array}{c|c}
\mathcal{A}_1     &  \mathcal{B}_1\\
\hline
\mathcal{C}_1     &  \mathcal{D}_1
\end{array}\right]$
and 
$\left[\begin{array}{c|c}
\mathcal{A}_2     &  \mathcal{B}_2\\
\hline
\mathcal{C}_2     &  \mathcal{D}_2
\end{array}\right]$
can be described by the LTI system
$\left[\begin{array}{c|c}
\begin{array}{cc}
\mathcal{A}_1   & \mathcal{B}_1 \mathcal{C}_2  \\
\mathbf{0}      & \mathcal{A}_2
\end{array}     
& 
\begin{array}{c}
\mathcal{B}_1 \mathcal{D}_2   \\
\mathcal{B}_2     
\end{array}\\
\hline
\begin{array}{cc}
\mathcal{C}_1      & \mathcal{D}_1\mathcal{C}_2
\end{array}     
& 
\mathcal{D}_1\mathcal{D}_2
\end{array}\right]$.
The notation for LPV systems is analogous to the LTI case with the dependence on the scheduling parameter explicitly shown.


\section{Problem Setup}\label{sec:problem_setup}
We consider source-seeking scenarios where one or more vehicles moving in $\mathbb{R}^d$ space (typically $d=\{1,2,3\}$) are embedded in an underlying differentiable scalar field $f: \mathbb{R}^d \xrightarrow[]{}\mathbb{R}$ with some convexity properties (made precise later).
We assume that each vehicle has access to the gradient $\nabla f$ evaluated at the vehicle's position.
Since this paper deals primarily with performance analysis, we assume that a local tracking controller has been designed and the closed-loop dynamics of the vehicle with reference input denoted by $r(t)\in \mathbb{R}^{2d}$ with $r(t)=[r_{\textnormal{pos}}^T(t) \quad r_{\textnormal{vel}}^T(t)]^T$, state $x(t)\in \mathbb{R}^n$ and position output $y(t)\in \mathbb{R}^d$ can be represented by 
\begin{equation}\label{eq:vehicle_dyn}
    \begin{split}
        \Dot{x}(t)&=Ax(t) + Br(t), \quad\quad x(0)=x_0,\\
        y(t)&=Cx(t).
    \end{split}
\end{equation}
We propose to augment these closed-loop dynamics by the second order dynamics,
\begin{equation}\label{eq:vir_vehicle_dyn}
    \begin{split}
        \Dot{r}_{\textnormal{pos}}(t)&=r_{\textnormal{vel}}(t),\\
        \Dot{r}_{\textnormal{vel}}(t)&=-k_d \cdot r_{\textnormal{vel}}(t) - k_p \cdot u(t),
    \end{split}
\end{equation}
where, $u(t) \in \mathbb{R}^d$ denotes external input.
\begin{figure}[ht]
	\tikzstyle{block} = [draw, rectangle, 
    minimum height=3em, minimum width=6em]
\tikzstyle{sum} = [draw,  circle, node distance=1cm]
\tikzstyle{input} = [coordinate]
\tikzstyle{output} = [coordinate]
\tikzstyle{pinstyle} = [pin edge={to-,thin,black}]

%

\begin{tikzpicture}[auto, node distance=2cm,>=latex']
\node [input, name=input] {};
\node [block, right of=input,node distance=3cm](Flockcontrol) {$\begin{array}{cc}
   \Dot{r}_{\textnormal{pos}} =& r_{\textnormal{vel}} \\
   \Dot{r}_{\textnormal{vel}} =&-k_d r_{\textnormal{vel}}-k_p u
\end{array}$};
\node [block, right of=Flockcontrol,node distance=4cm,minimum height=2.5em,minimum width=3em](velocityloop) {$\left[\begin{array}{c|c}
A     &  B\\
\hline
C     &  \mathbf{0}
\end{array}\right]$};
\node [output, right of=velocityloop,node distance=1.4cm](output) {};
\node [block, below of=Flockcontrol,node distance=2cm](delta) {$u=\nabla f(y)$};
\node[text=red, above left= 5mm and 0.5mm of Flockcontrol] (veh) {$G$};
\draw[red, dashed] (veh.east)-|([xshift=-3.8mm]output.west)|-([yshift=-3mm]Flockcontrol.south)-|([xshift=2mm]input.west)|-(veh.west);
\draw [draw,->] (input) -- node {$u$} (Flockcontrol);
\draw [draw,->] (Flockcontrol) -- node {$\begin{bmatrix}r_{\textnormal{pos}}\\r_{\textnormal{vel}}\end{bmatrix}$} (velocityloop);
\draw [draw,->] (velocityloop) -- node[name=outmid] {$y$} (output);
\draw [->] (outmid) |- (delta);
\draw [-] (delta) -| (input);
\end{tikzpicture}
	\caption{Control architecture}
	\label{fig:control_architecture}	
\end{figure}
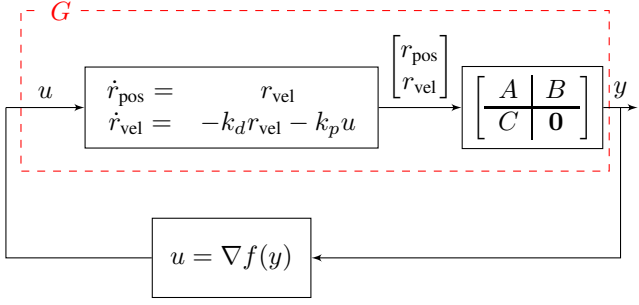
As shown in Fig. \ref{fig:control_architecture}, we consider a control architecture with $u(t)=\nabla f(y(t))$.
The intuition behind the control architecture is to see that if the tracking controller is doing a good job ($y(t)\approx r_{\textnormal{pos}}(t)$), we end up with a second order dissipative Hamiltonian system, the stability of which can be analyzed as in \cite{Datar.51220205152020} and \cite{Attallah.2020}. 
The overall vehicle dynamics from $u$ to $y$ denoted by $G$, along with the feedback law $u(t)=\nabla f(y(t))$,  can be represented by 
\begin{equation} \label{eq:sys_dyn_G}
    \begin{split}
        \Dot{\eta}(t)&=A_G\eta(t) + B_G u(t), \quad \quad \eta(0)=\eta_0,\\
        y(t)&=C_G \eta(t), \\
        u(t)&=\nabla f(y(t)),
    \end{split}
\end{equation}
where $\eta\in \mathbb{R}^{n_\eta}$ denotes the concatenated state vector.

We now make the following assumption.

\textbf{Assumption 1:} Let $f \in \mathcal{S}(m,L)$ and let $y_*$ minimize $f$, i.e., $f(y)\geq f(y_*) \, \forall y \in \mathbb{R}^d$ and $\nabla f(y_*)=0$.
\begin{remark}
The assumption on continuous differentiabiliity can be relaxed by using sub-gradient maps, but the assumption on strong convexity is essential.
\end{remark}
Since we consider vehicle models, we assume that $A_G$ has an eigen value at 0 (integral action).
Observe that this is enforced by the control architecture of $G$ with the first block having an integrator in series with a stable block. 
Since this implies that $\exists \eta_*\neq 0, A_G\eta_*=0$, we get a feasible equilibrium $(\eta_*,u_*,y_*)$ for \eqref{eq:sys_dyn_G} satisfying
\begin{equation} \label{eq:sys_dyn_G_delta}
    \begin{split}
        0&=A_G\eta_* + B_G u_*=A_G\eta_* ,\\
        y_*&=C_G \eta_* ,\\
        u_*&=\nabla f(y_*)=0.
    \end{split}
\end{equation}
The dynamics in the deviation variables $\Tilde{\eta}(t)=\eta(t)-\eta_*$, $\Tilde{u}(t)=u(t)-u_*=u(t)$ and $\Tilde{y}(t)=y(t)-y_*$ can be represented by  
\begin{equation} \label{eq:sys_dyn_G_tilde}
    \begin{split}
        \Dot{\Tilde{\eta}}(t)&=A_G\Tilde{\eta}(t) + B_G \Tilde{u}(t), \quad \quad \Tilde{\eta}(0)=\eta_0-\eta_*,\\
        \Tilde{y}(t)&=C_G \Tilde{\eta}(t)
    \end{split}
\end{equation}
and 
\begin{equation} \label{eq:delta}
\Tilde{u}(t)=\nabla f(\Tilde{y}(t)+y_*).
\end{equation}

\subsection{Flocking protocols under quadratic fields}
We now consider a group of $N$ identical LTI agents located in the scalar field with interactions among agents captured by flocking dynamics \cite{OlfatiSaber.2006} and along with a gradient-based forcing term driving them towards the source.
We use the notation and terminology from \cite{OlfatiSaber.2006}, \cite{Datar.51220205152020} and \cite{Attallah.2020} and refer the reader to these sources for details.
Let $\eta_i(t), u_i(t), y_i(t)$ denote the state, input and output of the $i^{\textnormal{th}}$ agent respectively.
We stack all the states of the agents to form a long state vector as $\hat{\eta}(t)=[ \eta_1(t)^T \cdots  \eta_N(t)^T ]^T$ and similarly define $\hat{u}(t)$ and $\hat{y}(t)$.
Let $\hat{f}:\mathbb{R}^{Nd} \rightarrow \mathbb{R}$ be defined as $\sum_{i=1}^N f(y_i)$ which results to $\nabla \hat{f}(\hat{y})=[ \nabla f(y_1)^T \cdots  \nabla f(y_1)^T ]^T$.
Let $\mathcal{V}:\mathbb{R}^{Nd} \rightarrow \mathbb{R}$ denote the interaction potential among different agents and let $\mathcal{L}$ denote the graph-Laplacian matrix representing the communication topology between agents.
Since the interaction between agents is defined to be such that there is force balance in every interacting pair of agents, it can be shown \cite{Datar.51220205152020} that,
\begin{equation}\label{eq:force_balance_flocking}
    \begin{split}
        (\mathbf{1}_N^T \otimes I_{d}) \cdot \nabla \mathcal{V}(\hat{y}) &= \mathbf{0} \quad \forall \hat{y} \in \mathbb{R}^{Nd},\\
        (\mathbf{1}_N^T \otimes I_{d}) (\mathcal{L}\otimes I_d)=(\mathbf{1}_N^T \mathcal{L})\otimes I_d &= \mathbf{0}.
    \end{split}
\end{equation}
The overall flocking dynamics with gradient based forcing terms could be written as
\begin{equation}
    \begin{split}
        \Dot{\hat{\eta}}(t)&=(I_N \otimes A_G)\hat{\eta}(t) + (I_N \otimes B_G)\hat{u}(t),\\
        \hat{y}(t)&=(I_N \otimes C_G)\hat{\eta}(t) ,\\
        \hat{u}(t)&=\nabla \hat{f}(\hat{y}(t))-\nabla \mathcal{V}(\hat{y}(t))-(\mathcal{L}\otimes I_d)\Dot{\hat{y}}(t).
    \end{split}
\end{equation}
Let $\eta_c=\frac{1}{N}(\mathbf{1}_N^T \otimes I_{n_\eta})\hat{\eta}$, $y_c=\frac{1}{N}(\mathbf{1}_N^T \otimes I_d)\hat{y}$ and $u_c=\frac{1}{N}(\mathbf{1}_N^T \otimes I_d)\hat{u}$ be the COM state, input and output, respectively.

If $f$ is a quadratic field of the form $f(y)=y^T Q y + c^Ty + d$, we get, $\nabla f (y)=2 Qy + c$ which is linear.
Therefore, $\nabla \hat{f}(\hat{y})=(I_N \otimes Q)\hat{y}+ \mathbf{1}_N \otimes c$ and $\frac{1}{N}(\mathbf{1}_N^T \otimes I_{d})\nabla \hat{f}(\hat{y}(t))=\nabla f (y_c)$.
Using this fact and as well as \eqref{eq:force_balance_flocking}, we can derive the reduced center of mass (COM) dynamics as
\begin{equation}\label{eq:COM_dynamics}
    \begin{split}
        \Dot{\eta_c}(t)
        &=A_G \eta_c + B_G u_c,\\
        y_c(t)
        &= C_G \eta_c ,\\
        u_c(t)&=\nabla f (y_c).
    \end{split}
\end{equation}
This has the same form as \eqref{eq:sys_dyn_G} and hence, the analysis of a single agent can be used to analyze the performance of multiple agents flocking under quadratic scalar fields.
\begin{remark}
It is important to note that structural dynamics (disagreement dynamics) has to be analyzed separately such as in \cite{OlfatiSaber.2006} to ensure the the stability of the overall flocking system. We focus in this paper on the macroscopic performance of the group for source-seeking, which as shown above reduces to the analysis of a single agent.
\end{remark}
\section{Theory}\label{sec:theory}
\begin{figure}[t]
	\tikzstyle{block} = [draw, rectangle, 
    minimum height=3em, minimum width=6em]
\tikzstyle{sum} = [draw,  circle, node distance=1cm]
\tikzstyle{input} = [coordinate]
\tikzstyle{output} = [coordinate]
\tikzstyle{pinstyle} = [pin edge={to-,thin,black}]

\begin{tikzpicture}[auto, node distance=2cm,>=latex']
    \node [block,xshift=-1cm] (controller) 
    {$\begin{bmatrix}
    -mI_d  & I_d \\
    LI_d    & -I_d
    \end{bmatrix}$};
    \node [block, right of=controller,xshift=2cm,yshift=0.6cm] (filter1) {$\left[e^{-2\alpha t}h(-t)\right]*p(t)$};
    \node [block, right of=controller,xshift=2cm,yshift=-0.6cm,minimum width=3.15cm] (filter2) {$\left[e^{-2\alpha t}h(t)\right]*q(t)$};
    \node [input, left=of controller.170,xshift=1cm] (input1) {};
    \node [input, left=of controller.190,xshift=1cm] (input2) {};
    \node [input, left=of filter1.west,xshift=1.5cm] (p1) {};
    \node [input, left=of filter2.west,xshift=1.5cm] (q1) {};
    \node [output, right of=filter1] (output1) {};
    \node [output, right of=filter2] (output2) {};

    \draw[->] (input1) -- node [name=ytilde] {$\Tilde{y}$}(controller.170);
    \draw[->] (input2) -- node [name=utilde] {$\Tilde{u}$}(controller.190);
    \draw[-] (controller.10) -| node [name=p] {}(p1);
    \draw[-] (controller.350) -| node [name=q] {}(q1);
    \draw[->] (p1) -- node [name=p] {$p$}(filter1.west);
    \draw[->] (q1) -- node [name=q] {$q$}(filter2.west);
    \draw[->] (filter1) -- node [name=w2tilde] {$w_2$}(output1);
    \draw[->] (filter2) -- node [name=w1tilde] {$w_1$}(output2);

\end{tikzpicture}
	\caption{Signal definitions}
	\label{fig:signal_definitions}	
\end{figure}
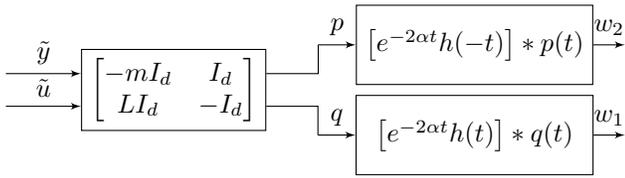

The fundamental idea in the framework of IQCs is to find properties (inequalities) satisfied by the input and output signals $\Tilde{u}$ and $\Tilde{y}$ of the gradient map $\Tilde{u}(t)=\nabla f(\Tilde{y}(t)+y_*)$.
It is convenient to derive properties of some new signals (defined below) related to the input and output signals $\Tilde{u}$ and $\Tilde{y}$. 
For the constants $m,L$ from Assumption 1, for any $\Tilde{u},\Tilde{y} \in \mathcal{L}_{2e}[0,\infty)$, define for $t\in [0,\infty)$,
\begin{equation}\label{eq:signal_defns}
    \begin{split}
        p(t)&=\Tilde{u}(t)-m \Tilde{y}(t),\\
        q(t)&=L \Tilde{y}(t) - \Tilde{u}(t).
    \end{split}
\end{equation}
With any $h \in \mathcal{L}_1(-\infty,\infty)$ satisfying
\begin{equation} \label{eq:zf_impulse_resp_cond}
    \begin{aligned}
        h(s) &\geq 0 \quad \forall s \in \mathbb{R} \textnormal{ and }\\
        \int_{-\infty}^{\infty}h(s)ds &\leq H,
    \end{aligned}
\end{equation}
let
\begin{equation} \label{eq:signal_defns_w12}
    \begin{split}
        w_1(t)&=\int_0^t e^{-2\alpha (t-\tau)}h(t-\tau) q(\tau) d\tau ,\\
        w_2(t)&=\int_0^t e^{-2\alpha (t-\tau)}h(-(t-\tau)) p(\tau) d\tau .
    \end{split}
\end{equation}
These signals have been shown in Fig. \ref{fig:signal_definitions} in the form of a block diagram.
We first present a key technical result which is covered by Lemma 3 from \cite{Freeman.62720186292018} where the result is presented in a very general setting of Bochner spaces.  
Moreover, since the proof of Lemma 3 from \cite{Freeman.62720186292018} is unavailable, we present a self-contained proof in the appendix adapted to our specific case building on ideas presented in a sidebar in \cite{Scherer.16052021}.
\begin{lemma}\label{theom:lemma_pq_ZF}
Let $\alpha \geq 0$ be fixed and let $\beta(\tau)=\textnormal{min}\{1,e^{-2\alpha \tau}\}$ for $\tau \in \mathbb{R}$. Then for all $\Tilde{u},\Tilde{y} \in \mathcal{L}_{2e}[0,\infty)$ that satisfy \eqref{eq:delta} with $f$ satisfying Assumption 1, the signals $p$ and $q$ as defined in \eqref{eq:signal_defns} satisfy, $\forall \tau \in \mathbb{R}, \forall T \geq 0$, 
\begin{equation}\label{eq:lemma_ineq}
    \int_0^T e^{2\alpha t} p(t)^T(q(t)-\beta(\tau)q_T(t-\tau))d\tau \geq 0,
\end{equation}
where, 
$q_T$ denotes the extension defined in \eqref{eq:extension_signal}.
\end{lemma}
\begin{proof}
See Appendix.
\end{proof}
Since \eqref{eq:lemma_ineq} holds for all $\tau \in \mathbb{R}$, we can conically combine the inequality \eqref{eq:lemma_ineq} by multiplying it  by a positive function $h(\tau)$ satisfying \eqref{eq:zf_impulse_resp_cond} and integrating out $\tau$. 
This leads to the next result.
\begin{theorem}\label{theom:theorem_pqw_filtered_ZF}
Let $h$ satisfying \eqref{eq:zf_impulse_resp_cond} be fixed and let $\alpha \geq 0$. For any $\Tilde{u}, \Tilde{y} \in \mathcal{L}_{2e}[0,\infty)$ that satisfy \eqref{eq:delta} with $f$ satisfying Assumption 1, the signals defined in \eqref{eq:signal_defns}, \eqref{eq:signal_defns_w12} satisfy
\begin{equation}
    \int_0^T e^{2\alpha t}( H p(t)^T q(t) - p(t)^T w_1(t) - q(t)^T w_2(t)) dt \geq 0,
\end{equation}
$\forall T \geq 0$.
\end{theorem}
\begin{proof}
See Appendix
\end{proof}
We now parameterize $h$ by proceeding along the lines of \cite{J.Veenman.2013}.
Let $A_{\nu}$ and $B_{\nu}$ be defined as
\begin{equation}
    A_{\nu}=
    \begin{bmatrix} 
    \lambda     & 0         &\dots     & 0 \\
    1           &\lambda    & \ddots   &0  \\
    0           & \ddots    & \ddots   &0  \\
    \vdots      & 0         & 1        &\lambda  
    \end{bmatrix},
    B_{\nu}=
    \begin{bmatrix} 
    1 \\
    0\\
    \vdots\\
    0 
    \end{bmatrix}.
\end{equation}
As pointed out in \cite{J.Veenman.2013}, it is not clear on how to choose $\lambda$ in the above parameterization.
For all our numerical experiments, we fix $\lambda=-1$ which gives reasonable results and further improvement in the results could be obtained by performing a line-search as suggested in \cite{J.Veenman.2013}.
Let $Q_{\nu}(t)=e^{A_{\nu}t}B_{\nu}=e^{\lambda t}R_{\nu} \left[ 1 \quad t \quad \hdots \quad t^{\nu-1}\right]^T$ where
$ R_{\nu}=\textnormal{diag}(\frac{1}{0!},\frac{1}{1!},\cdots,\frac{1}{(\nu-1)!})$.
Let $A_{\nu}^{\alpha}=A_{\nu}-2\alpha I$ and $Q_{\nu}^{\alpha}(t)=e^{-2\alpha t}Q_{\nu}(t)$.
Let $\Tilde{\psi}_{\nu}(s)=\left[ 1 \quad \frac{s}{(s-\lambda)^{\nu -1}} \quad \hdots \quad \frac{s^{\nu -1}}{(s-\lambda)^{\nu -1}}\right]^T$ and let $\left[
\begin{array}{c|c}
\Tilde{A}_{\nu} & \Tilde{B}_{\nu} \\
\hline
\Tilde{C}_{\nu} & \Tilde{D}_{\nu}
\end{array}
\right]$ be a state-space realization of $\Tilde{\psi}_{\nu}$.
Let the state space realization of $\Psi$ be as described in the appendix and shown in Fig. \ref{fig:psi_strucure}.
\begin{figure}[t]
	\tikzstyle{block} = [draw, rectangle, 
    minimum height=3em, minimum width=6em]
\tikzstyle{sum} = [draw,  circle, node distance=1cm]
\tikzstyle{input} = [coordinate]
\tikzstyle{output} = [coordinate]
\tikzstyle{pinstyle} = [pin edge={to-,thin,black}]

\begin{tikzpicture}[auto, node distance=2cm,>=latex']
    \node [block,xshift=-1cm] (controller) 
    {$\begin{bmatrix}
    -mI_d  & I_d \\
    LI_d    & -I_d
    \end{bmatrix}$};
    \node [block, right of=controller,xshift=2cm,yshift=0.7cm] (filter1) {$\left[e^{A_{\nu}^{\alpha}t}B_{\nu}\right]*p(t)$};
    \node [block, right of=controller,xshift=2cm,yshift=-0.7cm,minimum width=2.8cm] (filter2) {$\left[e^{A_{\nu}^{\alpha}t}B_{\nu}\right]*q(t)$};
    \node [input, left=of controller.170,xshift=1cm] (input1) {};
    \node [input, left=of controller.190,xshift=1cm] (input2) {};
    \node [input, left=of filter1.west,xshift=1.5cm] (p1) {};
    \node [input, left=of filter2.west,xshift=1.5cm] (q1) {};
    \node [output, right of=filter1] (output1) {};
    \node [output, right of=filter2] (output2) {};
    \node [draw,rectangle, right=of controller,xshift=2.5cm,fill=black,minimum height=10em, minimum width=0.05mm] (mux) {};
    \node [input, left=of mux.95,xshift=-1.5cm] (p11) {};
    \node [input, left=of mux.180,xshift=-1.5cm] (q11) {};
    \node [output, right=of mux,xshift=-1.5cm] (z) {};
    \node[text=red, above left= 11mm and 0.5mm of controller] (Psi) {$\Psi$};
    \draw[red, dashed] (Psi.east)-|([xshift=-3.8mm]z.west)|-([yshift=-14mm]controller.south)-|([xshift=2mm]input1.west)|-(Psi.west);

    \draw[->] (input1) -- node [name=ytilde] {$\Tilde{y}$}(controller.170);
    \draw[->] (input2) -- node [name=utilde] {$\Tilde{u}$}(controller.190);
    \draw[-] (controller.10) -| node [name=p] {}(p1);
    \draw[-] (controller.350) -| node [name=q] {}(q1);
    \draw[->] (p1) -- node [name=p] {$p$}(filter1.west);
    \draw[->] (q1) -- node [name=q] {$q$}(filter2.west);
    \draw[-] (p1) |- node [] {}(p11);
    \draw[-] (q1) |- node [] {}(q11);
    \draw[->] (p11) -- node [name=p11] {}(mux.95);
    \draw[->] (q11) -- node [name=q11] {}(mux.180);
    \draw[->] (filter1) -- node [name=w2tilde] {}(mux.100);
    \draw[->] (filter2) -- node [name=w1tilde] {}(mux.260);
    \draw[->] (mux) -- node [name=z] {$\Tilde{z}$}(z);

\end{tikzpicture}
	\caption{Structure of $\Psi$}
	\label{fig:psi_strucure}	
\end{figure}
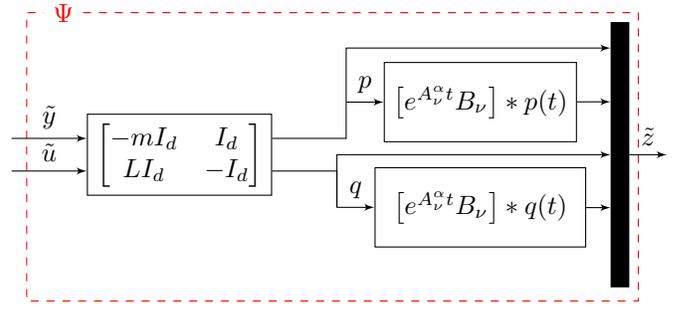
Before we present the next result, let us define for any $\Tilde{u}, \Tilde{y} \in \mathcal{L}_{2e}[0,\infty)$,
\begin{equation}\label{eq:signal_defn_z}
    \Tilde{z}(t)=\int_0^t C_{\Psi} e^{A_{\Psi}(t-\tau)}B_{\Psi}
\begin{bmatrix}
\Tilde{y}(\tau)\\
\Tilde{u}(\tau)
\end{bmatrix}d\tau
+ 
D_{\Psi}
\begin{bmatrix}
\Tilde{y}(t)\\
\Tilde{u}(t)
\end{bmatrix}.
\end{equation}
Consider the constraint in variables $H \in \mathbb{R},P_1\in \mathbb{R}^{1\times \nu},P_3\in \mathbb{R}^{1 \times \nu}$,
\begin{equation}\label{eq:L1_norm_constraint}
H+(P_1+P_3)A_{\nu}^{-1}B_{\nu} \geq 0,
\end{equation}
and consider the condition on $P_1,P_3$
\begin{equation}\label{eq:positivity}
\begin{split}
    \exists \mathcal{X}_1, \mathcal{X}_3 \in \mathbb{S}^{\nu-1}& \textnormal{such that for } i\in\{1,3\},\\
    (*)
    \begin{bmatrix}
    \mathbf{0}    & \mathcal{X}_i & \mathbf{0} \\
    \mathcal{X}_i & \mathbf{0}    & \mathbf{0} \\
    \mathbf{0} & \mathbf{0} & \textnormal{diag}(P_i)
    \end{bmatrix}
    &
    \begin{bmatrix}
    I                       & \mathbf{0} \\
    \Tilde{A}_{\nu}         & \Tilde{B}_{\nu} \\
    R_{\nu}\Tilde{C}_{\nu}  & R_{\nu}\Tilde{D}_{\nu} \\
    \end{bmatrix}
    >0.
\end{split}
\end{equation}
Let 
\begin{equation}
\mathbb{P}=\left\{
\begin{bmatrix}
\mathbf{0}  & \begin{bmatrix}
H  & -P_3 \\
-P_1^T    & \mathbf{0}
\end{bmatrix} \\
*    & \mathbf{0}
\end{bmatrix}:H, P_1, P_3\textnormal{ satisfy \eqref{eq:L1_norm_constraint},\eqref{eq:positivity}}\right\}.
\end{equation}

\begin{theorem}\label{theom:theorem_LMI_ZF}
For any $\Tilde{u},\Tilde{y} \in \mathcal{L}_{2e}[0,\infty)$ that satisfy \eqref{eq:delta} with $f$ satisfying Assumption 1, the signal $\Tilde{z}$ as defined in \eqref{eq:signal_defn_z} satisfies the $\alpha-$IQCs
\begin{equation}
    \int_0^T e^{2\alpha t}\Tilde{z}^T(t) (P \otimes I_d) \Tilde{z}(t) dt \geq 0 \quad \forall P \in \mathbb{P},\forall T \geq 0.
\end{equation}
\end{theorem}
\begin{proof}
See Appendix
\end{proof}
\begin{remark}\label{rem:CC_ZF_cusal}
The function $h$ is usually referred to as the multiplier and by enforcing $P_1=\mathbf{0}$ ($P_3=\mathbf{0}$), we restrict the search space to causal (anti-causal) ZF multipliers.
By enforcing $P_1=\mathbf{0}$ and $P_3=\mathbf{0}$, we specialize to the case of static multipliers which corresponds to the well-known circle-criterion (CC) \cite{Scherer.16052021}.
The conservatism between these specialization is investigated in section \ref{sec:numerical_results}.
\end{remark}
\begin{remark}
Extension of these results to cases when the map from $\Tilde{y}$ to $\Tilde{u}$ is additionally known to be odd is possible along the same lines but is not pursued in this paper. 
\end{remark}

We now present the final analysis result which leads to a computational procedure for obtaining convergence rate estimates.
Let $\left[\begin{array}{c|c}
\mathcal{A}     &  \mathcal{B}\\
\hline
\mathcal{C}     &  \mathcal{D}
\end{array}\right]$ be a state-space representation of the LTI system formed by the series interconnection of $\Psi$ and $\begin{bmatrix}
G\\
I
\end{bmatrix}$.
The following theorem from \cite{Hu.2016} gives the final performance analysis condition.
\begin{theorem}[\cite{Hu.2016}]\label{theom:theorem_perf_analysis}
If $\exists \mathcal{X}>0,P \in \mathbb{P}$ such that
\begin{equation}\label{eq:perf_LMI}
\begin{bmatrix}
\mathcal{A}^T\mathcal{X}+\mathcal{X}\mathcal{A}+2\alpha \mathcal{X}  & \mathcal{X}\mathcal{B} \\
\mathcal{B}^T\mathcal{X}    & \mathbf{0}
\end{bmatrix}
+
\begin{bmatrix}
\mathcal{C}^T \\
\mathcal{D}^T 
\end{bmatrix}
( P\otimes I_d)
\begin{bmatrix}
\mathcal{C} & \mathcal{D} 
\end{bmatrix}
\leq 
0,
\end{equation}
then, under dynamics \eqref{eq:sys_dyn_G} with $f$ satisfying Assumption 1, $y$ converges exponentially to $y_*$ with rate $\alpha$, i.e., $\exists \kappa \geq 0$ such that  
$||\Tilde{y}(t)||\leq \kappa e^{-\alpha t}$ holds for all $t\geq 0$.
\end{theorem}
\begin{proof}
See \cite{Hu.2016}.
\end{proof}
\begin{remark}
We point out that \eqref{eq:perf_LMI} is not linear in $\alpha$, $\mathcal{X}$ and $P$ due to the product $\alpha \mathcal{X}$.
It falls into the class of quasi-convex optimization problems which can be solved efficiently.
We perform a bisection over $\alpha$ as suggested in \cite{Lessard.2016}.
\end{remark}

\section{Extension to LPV systems}\label{sec:lpv_ext}
Extensions of the results obtained in the previous section to LPV/ uncertain systems is straightforward and we demonstrate one such extension next.
Instead of the LTI system $G$, let $G(\rho)$ denote an LPV system with $n_{\rho}$ scheduling parameters \cite{shamma1992linear}, where, for a compact set $\mathcal{P}\subset \mathbb{R}^{n_{\rho}} $, the function $\rho:[0,\infty) \rightarrow \mathcal{P}$ captures the time-dependence of the model parameters.
The dynamics of $G(\rho)$ can be represented by 
\begin{equation} \label{eq:sys_dyn_G_lpv}
    \begin{split}
        \Dot{\eta}(t)&=A_G(\rho(t))\eta(t) + B_G(\rho(t)) u(t), \quad \quad \eta(0)=\eta_0,\\
        y(t)&=C_G(\rho(t)) \eta(t), \\
        u(t)&=\nabla f(y(t)),
    \end{split}
\end{equation}
where $\eta\in \mathbb{R}^{n_\eta}$ is the state vector and $\rho:[0,\infty) \rightarrow \mathcal{P}$ is an arbitrary scheduling trajectory.
\begin{remark}
If the rate of parameter variation $\Dot{\rho}$ is bounded and this bound is known, we could include this information and consider parameter dependent Lyapunov functions (See \cite{.2000scherer} for details.), but we will not treat this case here.
\end{remark}

As for the linear case, we assume integral action in $G(\rho)$, that is, $\exists \eta_*\neq 0, A_G(\Bar{\rho})\eta_*=0 \quad \forall \Bar{\rho} \in \mathcal{P}$.
Note again that this is implied by the control architecture.
Following the same steps as in the LTI case, we can write the dynamics in the deviation variables analogous to \eqref{eq:sys_dyn_G_delta} and \eqref{eq:delta}.
Let the series interconnection of the LTI system $\Psi$ and the LPV system $G(\rho)$ be denoted by $\left[\begin{array}{c|c}
\mathcal{A}(\rho)     &  \mathcal{B}(\rho)\\
\hline
\mathcal{C}(\rho)     &  \mathcal{D}(\rho)
\end{array}\right] = \Psi \cdot  \begin{bmatrix}
G(\rho)\\
I
\end{bmatrix}$.
The following theorem gives a sufficient condition for performance analysis.
\begin{theorem}\label{theom:theorem_perf_analysis_lpv}
If $\exists \mathcal{X}>0,P \in \mathbb{P}$ such that, for any $\bar{\rho} \in \mathcal{P}$,
\begin{equation}\label{eq:perf_LMI_lpv}
\begin{split}
\begin{bmatrix}
\mathcal{A}(\bar{\rho})^T\mathcal{X}+\mathcal{X}\mathcal{A}(\bar{\rho})+2\alpha \mathcal{X}  & \mathcal{X}\mathcal{B}(\bar{\rho}) \\
\mathcal{B}(\bar{\rho})^T\mathcal{X}    & \mathbf{0}
\end{bmatrix}
+\\
\begin{bmatrix}
\mathcal{C}(\bar{\rho})^T \\
\mathcal{D}(\bar{\rho})^T 
\end{bmatrix}
( P\otimes I_d)
\begin{bmatrix}
\mathcal{C}(\bar{\rho}) & \mathcal{D}(\bar{\rho}) 
\end{bmatrix}
\leq 
0,   
\end{split}
\end{equation}
then, under dynamics \eqref{eq:sys_dyn_G_lpv} with $f$ satisfying Assumption 1, $y$ converges exponentially to $y_*$ with rate $\alpha$, i.e., $\exists \kappa \geq 0$ such that  
$||\Tilde{y}(t)||\leq \kappa e^{-\alpha t}$ holds for all $t\geq 0$.
\end{theorem}
\begin{proof}
See Appendix.
\end{proof}
\begin{remark}
We note that the equivalence between the center of mass dynamics and the dynamics of a single agent illustrated for homogeneous LTI agent flocking dynamics under quadratic fields in section \ref{sec:problem_setup} extends to LPV systems only if all agents are homogeneously scheduled, i.e., $\rho_i(t)=\rho_j(t) \quad \forall t \in [0,\infty),\quad \forall i,j \in \{1,2,\cdots, N\}$.
\end{remark}
\section{Numerical Results}\label{sec:numerical_results}
The code used for generating results in this section is availabe at \cite{ACC2022code}.
\subsection{Quadrotor}\label{eg:quadrotor}
\begin{figure}[t]
%
%
\definecolor{mycolor1}{rgb}{0.00000,1.00000,1.00000}%
\begin{tikzpicture}

\begin{axis}[%
width=2.521in,
height=2.366in,
at={(0.758in,0.481in)},
scale only axis,
xmin=1,
xmax=9,
xlabel style={font=\color{white!15!black}},
xlabel={L},
ymin=0,
ymax=0.5,
ylabel style={font=\color{white!15!black}},
ylabel={$\alpha$},
axis background/.style={fill=white},
title style={font=\bfseries},
title={},
legend style={legend cell align=left, align=left, draw=white!15!black}
]
\addplot [color=red, line width=1.0pt, draw=none, mark=o, mark options={solid, red}]
  table[row sep=crcr]{%
1	0.4119873046875\\
1.2	0.3472900390625\\
1.4	0.3076171875\\
1.6	0.277099609375\\
1.8	0.2508544921875\\
2	0.2276611328125\\
2.2	0.2069091796875\\
2.4	0.1873779296875\\
2.6	0.1690673828125\\
2.8	0.152587890625\\
3	0.1361083984375\\
3.2	0.120849609375\\
3.4	0.106201171875\\
3.6	0.0921630859375\\
3.8	0.0787353515625\\
4	0.06591796875\\
4.2	0.0531005859375\\
4.4	0.0408935546875\\
4.6	0.029296875\\
4.8	0.0177001953125\\
5	0.006103515625\\
5.2	-1\\
5.4	-1\\
5.6	-1\\
5.8	-1\\
6	-1\\
6.2	-1\\
6.4	-1\\
6.6	-1\\
6.8	-1\\
7	-1\\
7.2	-1\\
7.4	-1\\
7.6	-1\\
7.8	-1\\
8	-1\\
8.2	-1\\
8.4	-1\\
8.6	-1\\
8.8	-1\\
9	-1\\
};
\addlegendentry{CC($P_1=0,P_3=0$)}

\addplot [color=green, dashed, line width=1.0pt]
  table[row sep=crcr]{%
1	0.4119873046875\\
1.2	0.3643798828125\\
1.4	0.33935546875\\
1.6	0.321044921875\\
1.8	0.306396484375\\
2	0.2935791015625\\
2.2	0.281982421875\\
2.4	0.27099609375\\
2.6	0.260009765625\\
2.8	0.2490234375\\
3	0.2386474609375\\
3.2	0.2276611328125\\
3.4	0.21728515625\\
3.6	0.2069091796875\\
3.8	0.196533203125\\
4	0.1861572265625\\
4.2	0.17578125\\
4.4	0.1654052734375\\
4.6	0.1556396484375\\
4.8	0.145263671875\\
5	0.135498046875\\
5.2	0.125732421875\\
5.4	0.115966796875\\
5.6	0.1068115234375\\
5.8	0.0970458984375\\
6	0.087890625\\
6.2	0.0787353515625\\
6.4	0.069580078125\\
6.6	0.06103515625\\
6.8	0.0518798828125\\
7	0.0433349609375\\
7.2	0.0347900390625\\
7.4	0.0262451171875\\
7.6	0.0177001953125\\
7.8	0.009765625\\
8	0.0018310546875\\
8.2	-1\\
8.4	-1\\
8.6	-1\\
8.8	-1\\
9	-1\\
};
\addlegendentry{ZF causal ($P_1=0$)}

\addplot [color=blue, dashdotted, line width=1.0pt]
  table[row sep=crcr]{%
1	0.4119873046875\\
1.2	0.3472900390625\\
1.4	0.3076171875\\
1.6	0.277099609375\\
1.8	0.2508544921875\\
2	0.2276611328125\\
2.2	0.2069091796875\\
2.4	0.1873779296875\\
2.6	0.1690673828125\\
2.8	0.152587890625\\
3	0.1361083984375\\
3.2	0.120849609375\\
3.4	0.106201171875\\
3.6	0.0921630859375\\
3.8	0.0787353515625\\
4	0.06591796875\\
4.2	0.0531005859375\\
4.4	0.0408935546875\\
4.6	0.029296875\\
4.8	0.0177001953125\\
5	0.006103515625\\
5.2	-1\\
5.4	-1\\
5.6	-1\\
5.8	-1\\
6	-1\\
6.2	-1\\
6.4	-1\\
6.6	-1\\
6.8	-1\\
7	-1\\
7.2	-1\\
7.4	-1\\
7.6	-1\\
7.8	-1\\
8	-1\\
8.2	-1\\
8.4	-1\\
8.6	-1\\
8.8	-1\\
9	-1\\
};
\addlegendentry{ZF anti-causal ($P_3=0$)}

\addplot [color=mycolor1, line width=1.0pt]
  table[row sep=crcr]{%
1	0.4119873046875\\
1.2	0.3955078125\\
1.4	0.379638671875\\
1.6	0.3643798828125\\
1.8	0.34912109375\\
2	0.3338623046875\\
2.2	0.3192138671875\\
2.4	0.30517578125\\
2.6	0.2911376953125\\
2.8	0.2777099609375\\
3	0.2642822265625\\
3.2	0.25146484375\\
3.4	0.2386474609375\\
3.6	0.2264404296875\\
3.8	0.2142333984375\\
4	0.20263671875\\
4.2	0.1910400390625\\
4.4	0.179443359375\\
4.6	0.16845703125\\
4.8	0.157470703125\\
5	0.146484375\\
5.2	0.1361083984375\\
5.4	0.125732421875\\
5.6	0.1153564453125\\
5.8	0.1055908203125\\
6	0.0958251953125\\
6.2	0.0860595703125\\
6.4	0.0762939453125\\
6.6	0.067138671875\\
6.8	0.0579833984375\\
7	0.048828125\\
7.2	0.040283203125\\
7.4	0.0311279296875\\
7.6	0.0225830078125\\
7.8	0.0140380859375\\
8	0.006103515625\\
8.2	-1\\
8.4	-1\\
8.6	-1\\
8.8	-1\\
9	-1\\
};
\addlegendentry{ZF}

\addplot [color=black, line width=1.0pt, draw=none, mark=asterisk, mark options={solid, black}]
  table[row sep=crcr]{%
1	0.412331651029501\\
1.2	0.395974289263064\\
1.4	0.380000019452199\\
1.6	0.364400808966183\\
1.8	0.349167588131177\\
2	0.334290538830461\\
2.2	0.319759335099553\\
2.4	0.305563341010742\\
2.6	0.291691771302537\\
2.8	0.27813382002905\\
3	0.264878762120048\\
3.2	0.251916032253449\\
3.4	0.239235284914788\\
3.6	0.226826438996403\\
3.8	0.2146797097987\\
4	0.202785630850186\\
4.2	0.191135067568535\\
4.4	0.179719224441611\\
4.6	0.168529647113125\\
4.8	0.157558220508057\\
5	0.146797163923227\\
5.2	0.136239023833451\\
5.4	0.125876665018557\\
5.6	0.115703260496915\\
5.8	0.105712280652737\\
6	0.095897481864014\\
6.2	0.08625289487244\\
6.4	0.0767728130832965\\
6.6	0.067451780940349\\
6.8	0.0582845824858749\\
7	0.0492662301880655\\
7.2	0.0403919540956929\\
7.4	0.0316571913620181\\
7.6	0.0230575761658658\\
7.8	0.0145889300465953\\
8	0.00624725266099713\\
8.2	-0.00197128703663169\\
8.4	-0.0100703591950899\\
8.6	-0.0180534810554159\\
8.8	-0.0259240245392554\\
9	-0.0336852234492348\\
};
\addlegendentry{Example fields}

\end{axis}
\end{tikzpicture}%
	\caption{Convergence rate estimates for quadrotor dynamics provided by different multipliers(See Remark \ref{rem:CC_ZF_cusal}) for fields $f\in \mathcal{S}(1,L)$}
	\centering
	\label{fig:quadrotor_robustness}
\end{figure}
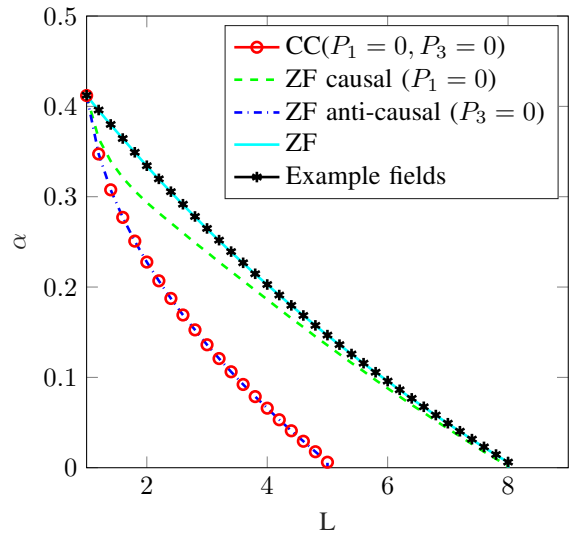
\begin{figure}[t]
%
%
\definecolor{mycolor1}{rgb}{0.00000,1.00000,1.00000}%
\begin{tikzpicture}

\begin{axis}[%
width=2.521in,
height=2.366in,
at={(0.758in,0.481in)},
scale only axis,
xmin=5,
xmax=30,
xlabel style={font=\color{white!15!black}},
xlabel={$\frac{k_d}{k_p}$},
ymin=0,
ymax=0.2,
ylabel style={font=\color{white!15!black}},
ylabel={$\alpha$},
axis background/.style={fill=white},
title style={font=\bfseries},
title={},
legend style={legend cell align=left, align=left, draw=white!15!black}
]
\addplot [color=red, line width=1.0pt, draw=none, mark=o, mark options={solid, red}]
  table[row sep=crcr]{%
1	-1\\
1.5	-1\\
2	-1\\
2.5	-1\\
3	-1\\
3.5	-1\\
4	-1\\
4.5	-1\\
5	-1\\
5.5	-1\\
6	-1\\
6.5	-1\\
7	-1\\
7.5	-1\\
8	-1\\
8.5	-1\\
9	-1\\
9.5	-1\\
10	-1\\
10.5	-1\\
11	0.006103515625\\
11.5	0.0140380859375\\
12	0.0213623046875\\
12.5	0.028076171875\\
13	0.0341796875\\
13.5	0.0390625\\
14	0.0439453125\\
16	0.05615234375\\
18	0.0592041015625\\
20	0.0543212890625\\
22	0.048828125\\
24	0.0445556640625\\
26	0.0408935546875\\
28	0.037841796875\\
30	0.0347900390625\\
};
\addlegendentry{CC($P_1=0,P_3=0$)}

\addplot [color=green, line width=1.0pt, draw=none, mark=asterisk, mark options={solid, green}]
  table[row sep=crcr]{%
1	-1\\
1.5	-1\\
2	-1\\
2.5	-1\\
3	-1\\
3.5	-1\\
4	-1\\
4.5	-1\\
5	-1\\
5.5	-1\\
6	-1\\
6.5	0.0244140625\\
7	0.0494384765625\\
7.5	0.072021484375\\
8	0.0921630859375\\
8.5	0.10986328125\\
9	0.125732421875\\
9.5	0.130615234375\\
10	0.1226806640625\\
10.5	0.1153564453125\\
11	0.108642578125\\
11.5	0.1025390625\\
12	0.09765625\\
12.5	0.0927734375\\
13	0.0885009765625\\
13.5	0.0848388671875\\
14	0.0811767578125\\
16	0.069580078125\\
18	0.06103515625\\
20	0.0543212890625\\
22	0.048828125\\
24	0.0445556640625\\
26	0.0408935546875\\
28	0.037841796875\\
30	0.0347900390625\\
};
\addlegendentry{ZF causal ($P_1=0$)}

\addplot [color=blue, dashdotted, line width=1.0pt]
  table[row sep=crcr]{%
1	-1\\
1.5	-1\\
2	-1\\
2.5	-1\\
3	-1\\
3.5	-1\\
4	-1\\
4.5	-1\\
5	-1\\
5.5	-1\\
6	-1\\
6.5	-1\\
7	-1\\
7.5	-1\\
8	-1\\
8.5	-1\\
9	-1\\
9.5	-1\\
10	-1\\
10.5	-1\\
11	0.006103515625\\
11.5	0.0140380859375\\
12	0.0213623046875\\
12.5	0.028076171875\\
13	0.0341796875\\
13.5	0.0390625\\
14	0.0439453125\\
16	0.05615234375\\
18	0.0592041015625\\
20	0.0543212890625\\
22	0.048828125\\
24	0.0445556640625\\
26	0.0408935546875\\
28	0.037841796875\\
30	0.0347900390625\\
};
\addlegendentry{ZF anti-causal ($P_3=0$)}

\addplot [color=mycolor1, line width=1.0pt]
  table[row sep=crcr]{%
1	-1\\
1.5	-1\\
2	-1\\
2.5	-1\\
3	-1\\
3.5	-1\\
4	-1\\
4.5	-1\\
5	-1\\
5.5	-1\\
6	0.001220703125\\
6.5	0.0311279296875\\
7	0.0579833984375\\
7.5	0.081787109375\\
8	0.103759765625\\
8.5	0.123291015625\\
9	0.140380859375\\
9.5	0.130615234375\\
10	0.1226806640625\\
10.5	0.1153564453125\\
11	0.108642578125\\
11.5	0.1025390625\\
12	0.09765625\\
12.5	0.0927734375\\
13	0.0885009765625\\
13.5	0.0848388671875\\
14	0.0811767578125\\
16	0.069580078125\\
18	0.06103515625\\
20	0.0543212890625\\
22	0.048828125\\
24	0.0445556640625\\
26	0.0408935546875\\
28	0.037841796875\\
30	0.0347900390625\\
};
\addlegendentry{ZF}

\end{axis}
\end{tikzpicture}%
	\caption{Performance for different $\frac{k_d}{k_p}$ for a fixed scalar field}
	\centering
	\label{fig:quadrotor_optimal_kd}
\end{figure}
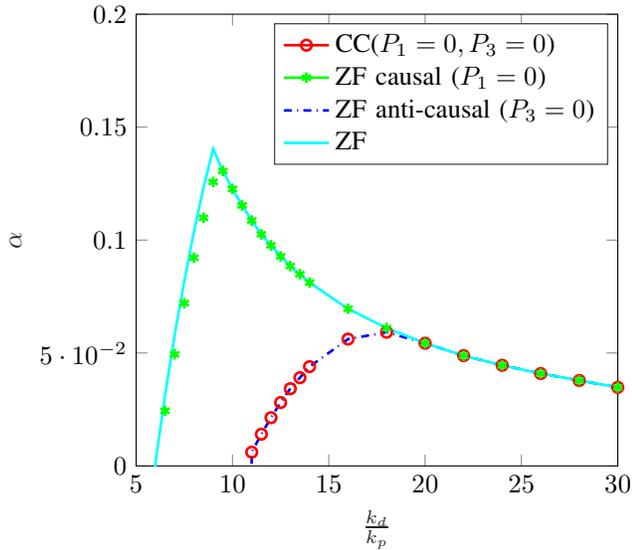
\begin{figure}[t]
	\input{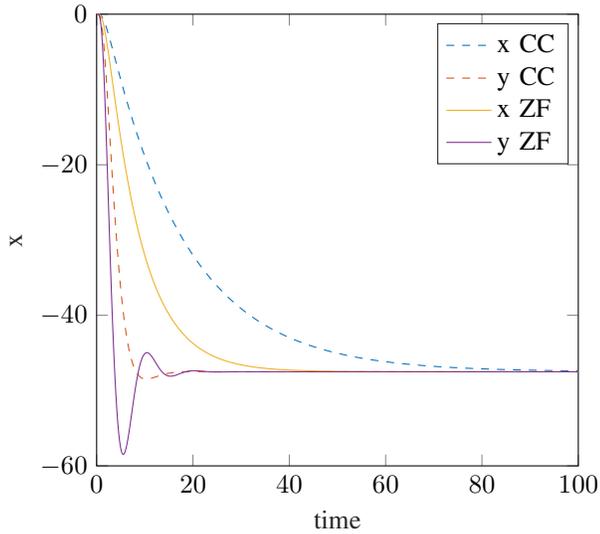}
	\caption{Sample trajectories of a quadrotor in a 2D-field with minimum at (-50,-50) with tuning based on CC (dashed lines) and tuning based on ZF (solid lines)}
	\centering
	\label{fig:quadrotor_trajectory}
\end{figure}

We consider a linearized quadrotor model and use an LQR (Linear-Quadratic-Regulator) based state-feedback controller tuned for zero steady-state error for step references.  
As discussed in section \ref{sec:problem_setup}, we let this closed-loop system be represented by state-space realization \eqref{eq:vehicle_dyn} and augment it with dynamics \eqref{eq:vir_vehicle_dyn} to obtain the system $G$. 

We address the following questions next to demonstrate the applicability of the theoretical results.
\begin{itemize}
    \item[1] How robust is the given controller with respect to different fields $f\in \mathcal{S}(m,L)$?
    \item[2] How do we design the gains $k_p$ and $k_d$ for the given closed-loop quadrotor system?
    \item[3] How conservative are the estimates of the convergence rates given by our analysis for static multipliers (circle criterion) and causal/anti-causal/non-causal ZF mutipliers (See Remark \ref{rem:CC_ZF_cusal}) ?
\end{itemize}
For fixed gains $k_p$ and $k_d$ and given closed-loop quadrotor dynamics, Fig. \ref{fig:quadrotor_robustness} shows the convergence rate estimates provided by different multipliers for fields $f\in \mathcal{S}(1,L)$ with increasing $L$.
Since increasing $L$ enlarges the set of allowable fields, i.e., $\mathcal{S}(1,L_1) \subset \mathcal{S}(1,L_2) \quad \forall L_1\leq L_2$, the estimates are non-increasing with increasing $L$. 
It can be seen that while we can certify stability with the cricle criterion for fields $f \in \mathcal{S}(1,5)$, the general non-causal ZF multipliers along with the ZF multipliers restricted to the causal case($P_3=0$) can certify stability for all fields $f \in \mathcal{L}(1,8)$.
Furthermore, for each $L$, we can find a field (quadratic) that achieves the convergence rate guaranteed by the analysis showing that, in this example, the estimates are tight.
The conservatism incurred by restricting the search to causal multipliers is minor in this example.
The stability analysis in \cite{Attallah.2020} uses manually constructed diagonal storage functions together with a small gain argument and for this example, gives the sufficient condition for stability to be $L<5$.
This interestingly coincides with the stability boundary given by the circle-criterion (static multipliers).
Since performance analysis was not included in \cite{Attallah.2020}, this example illustrates the extension of \cite{Attallah.2020} to a non-conservative performance analysis. 

The effect of varying the ratio of gains $\frac{k_d}{k_p}$ on the performance estimates for fixed allowable field set is shown in Fig. \ref{fig:quadrotor_optimal_kd}. 
It shows that the highest convergence rate of  $0.14$ can be achieved for $\frac{k_d}{k_p}=9$ and demonstrates a method for tuning the gains for optimal convergence rates.
We also observe that tuning the gains by using static multipliers (circle criterion) would lead to to a rather poor performance.
Sample trajectories of a quadrotor locating the source at (-50,-50) is shown in Fig. \ref{fig:quadrotor_trajectory}, where, $\frac{k_d}{k_p}$ are chosen optimal with respect to the circle criterion (dashed lines) and ZF (solid lines).
Although leading to a higher overshoot, the gains tuned with respect to ZF lead to faster convergence.

\subsection{Example showing the benefit of non-causal multipliers}
We now present an academic example that brings out the benefit of using general non-causal multipliers over causal multipliers.
Let $G(s)=5\frac{(s-1)}{s(s^2+s+25)}$ and consider fields $f\in \mathcal{S}(1,L)$.
The convergence rate estimates provided by different multipliers for increasing $L$ is shown in Fig. \ref{fig:robustness_non_causal_multiplier}.
It can be seen that while the circle criterion and causal ZF multipliers certify stability for fields $f\in \mathcal{S}(1,1.9)$, the anti-causal ZF multipliers can certify stability for fields $f\in \mathcal{S}(1,2)$ and the general non-causal ZF multipliers can certify stability for fields $f\in \mathcal{S}(1,2.4)$.
Furthermore, we can find example fields(quadratic) which coincide with the convergence rate estimates showing that these estimates are tight.
\begin{remark}
Since $f\in \mathcal{S}(1,L_1) \subset \mathcal{S}(1,L_2)$ for $L_1\leq L_2$, the estimates are non-increasing.
Note that the gap between the actual convergence rates for example fields and the estimates obtained from non-causal ZF multipliers for $L\in \{ 1.1, 1.2, \cdots, 1.6 \} $ does not imply conservatism since the example field used at $L=1$ is included for any larger $L$.
\end{remark}
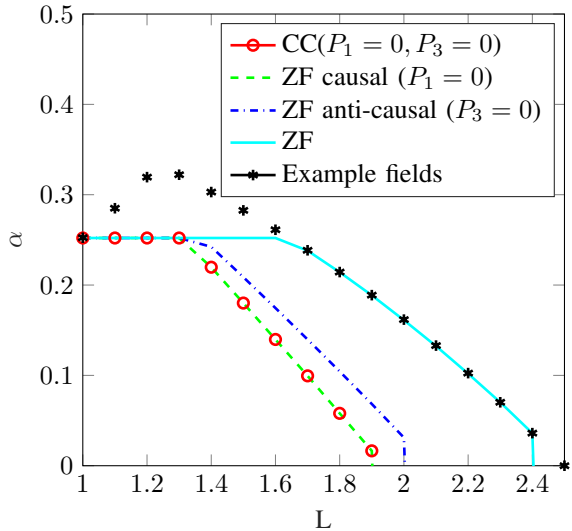
\begin{figure}[t]
%
%
\definecolor{mycolor1}{rgb}{0.00000,1.00000,1.00000}%
\begin{tikzpicture}

\begin{axis}[%
width=2.521in,
height=2.366in,
at={(0.758in,0.481in)},
scale only axis,
xmin=1,
xmax=2.5,
xlabel style={font=\color{white!15!black}},
xlabel={L},
ymin=0,
ymax=0.5,
ylabel style={font=\color{white!15!black}},
ylabel={$\alpha$},
axis background/.style={fill=white},
title style={font=\bfseries},
title={},
legend style={legend cell align=left, align=left, draw=white!15!black}
]
\addplot [color=red, line width=1.0pt, draw=none, mark=o, mark options={solid, red}]
  table[row sep=crcr]{%
1	0.2520751953125\\
1.1	0.2520751953125\\
1.2	0.2520751953125\\
1.3	0.2520751953125\\
1.4	0.2197265625\\
1.5	0.1800537109375\\
1.6	0.1397705078125\\
1.7	0.0994873046875\\
1.8	0.0579833984375\\
1.9	0.0164794921875\\
2	-1\\
2.1	-1\\
2.2	-1\\
2.3	-1\\
2.4	-1\\
2.5	-1\\
2.6	-1\\
2.7	-1\\
2.8	-1\\
2.9	-1\\
3	-1\\
};
\addlegendentry{CC($P_1=0,P_3=0$)}

\addplot [color=green, dashed, line width=1.0pt]
  table[row sep=crcr]{%
1	0.2520751953125\\
1.1	0.2520751953125\\
1.2	0.2520751953125\\
1.3	0.2520751953125\\
1.4	0.2197265625\\
1.5	0.1800537109375\\
1.6	0.1397705078125\\
1.7	0.0994873046875\\
1.8	0.0579833984375\\
1.9	0.0164794921875\\
2	-1\\
2.1	-1\\
2.2	-1\\
2.3	-1\\
2.4	-1\\
2.5	-1\\
2.6	-1\\
2.7	-1\\
2.8	-1\\
2.9	-1\\
3	-1\\
};
\addlegendentry{ZF causal ($P_1=0$)}

\addplot [color=blue, dashdotted, line width=1.0pt]
  table[row sep=crcr]{%
1	0.2520751953125\\
1.1	0.2520751953125\\
1.2	0.2520751953125\\
1.3	0.2520751953125\\
1.4	0.2423095703125\\
1.5	0.208740234375\\
1.6	0.174560546875\\
1.7	0.1397705078125\\
1.8	0.1043701171875\\
1.9	0.068359375\\
2	0.0311279296875\\
2.1	-1\\
2.2	-1\\
2.3	-1\\
2.4	-1\\
2.5	-1\\
2.6	-1\\
2.7	-1\\
2.8	-1\\
2.9	-1\\
3	-1\\
};
\addlegendentry{ZF anti-causal ($P_3=0$)}

\addplot [color=mycolor1, line width=1.0pt]
  table[row sep=crcr]{%
1	0.2520751953125\\
1.1	0.2520751953125\\
1.2	0.2520751953125\\
1.3	0.2520751953125\\
1.4	0.2520751953125\\
1.5	0.2520751953125\\
1.6	0.2520751953125\\
1.7	0.238037109375\\
1.8	0.2142333984375\\
1.9	0.1885986328125\\
2	0.1611328125\\
2.1	0.1324462890625\\
2.2	0.1019287109375\\
2.3	0.0701904296875\\
2.4	0.0360107421875\\
2.5	-1\\
2.6	-1\\
2.7	-1\\
2.8	-1\\
2.9	-1\\
3	-1\\
};
\addlegendentry{ZF}

\addplot [color=black, line width=1.0pt, draw=none, mark=asterisk, mark options={solid, black}]
  table[row sep=crcr]{%
1	0.252381023686084\\
1.1	0.285030034676752\\
1.2	0.319444587967699\\
1.3	0.322120595821607\\
1.4	0.30294136645159\\
1.5	0.282662461122388\\
1.6	0.261201184778797\\
1.7	0.238470096600161\\
1.8	0.214377762177947\\
1.9	0.188829987346274\\
2	0.161731674331652\\
2.1	0.132989449419385\\
2.2	0.102515199007407\\
2.3	0.0702306031940527\\
2.4	0.0360726564236493\\
2.5	-1.11022302462516e-16\\
2.6	-0.0380003387049316\\
2.7	-0.0779064749243448\\
2.8	-0.119656383851905\\
2.9	-0.163144853458106\\
3	-0.208223620537594\\
};
\addlegendentry{Example fields}

\end{axis}
\end{tikzpicture}%
	\caption{Robustness against different fields $f\in \mathcal{S}(1,L)$ for $G(s)=5\frac{(s-1)}{s(s^2+s+25)}$} 
	\centering
	\label{fig:robustness_non_causal_multiplier}
\end{figure}
\
\subsection{LPV generic vehicle model}
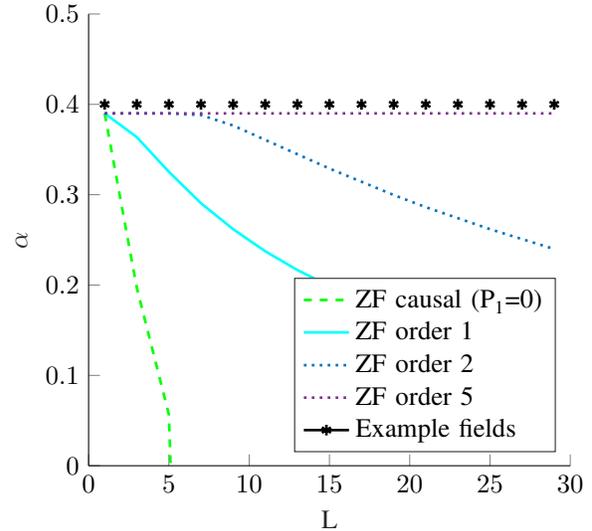
\begin{figure}[t]
%
%
\definecolor{mycolor1}{rgb}{0.00000,1.00000,1.00000}%
\definecolor{mycolor2}{rgb}{0.00000,0.44700,0.74100}%
\definecolor{mycolor3}{rgb}{0.85000,0.32500,0.09800}%
\definecolor{mycolor4}{rgb}{0.92900,0.69400,0.12500}%
\definecolor{mycolor5}{rgb}{0.49400,0.18400,0.55600}%
\begin{tikzpicture}

\begin{axis}[%
width=2.521in,
height=2.366in,
at={(0.758in,0.481in)},
scale only axis,
xmin=0,
xmax=30,
xlabel style={font=\color{white!15!black}},
xlabel={L},
ymin=0,
ymax=0.5,
ylabel style={font=\color{white!15!black}},
ylabel={$\alpha$},
axis background/.style={fill=white},
title style={font=\bfseries},
title={},
axis x line*=bottom,
axis y line*=left,
legend style={legend cell align=left, align=left, draw=white!15!black},
legend pos=south east
]
\addplot [color=green, dashed, line width=1.0pt]
  table[row sep=crcr]{%
1	0.3900146484375\\
3	0.1971435546875\\
5	0.0567626953125\\
7	-1\\
9	-1\\
11	-1\\
13	-1\\
15	-1\\
17	-1\\
19	-1\\
21	-1\\
23	-1\\
25	-1\\
27	-1\\
29	-1\\
};
\addlegendentry{$\text{ZF causal (P}_\text{1}\text{=0)}$}


\addplot [color=mycolor1, line width=1.0pt]
  table[row sep=crcr]{%
1	0.3900146484375\\
3	0.36376953125\\
5	0.3253173828125\\
7	0.29052734375\\
9	0.2618408203125\\
11	0.2374267578125\\
13	0.2166748046875\\
15	0.198974609375\\
17	0.1837158203125\\
19	0.1702880859375\\
21	0.15869140625\\
23	0.1483154296875\\
25	0.13916015625\\
27	0.130615234375\\
29	0.123291015625\\
};
\addlegendentry{ZF order 1}

\addplot [color=mycolor2, dotted, line width=1.0pt]
  table[row sep=crcr]{%
1	0.3900146484375\\
3	0.3900146484375\\
5	0.3900146484375\\
7	0.38818359375\\
9	0.3765869140625\\
11	0.3607177734375\\
13	0.3448486328125\\
15	0.3289794921875\\
17	0.3143310546875\\
19	0.2996826171875\\
21	0.2862548828125\\
23	0.2740478515625\\
25	0.2618408203125\\
27	0.2508544921875\\
29	0.2398681640625\\
};
\addlegendentry{ZF order 2}



\addplot [color=mycolor5, dotted, line width=1.0pt]
  table[row sep=crcr]{%
1	0.3900146484375\\
3	0.3900146484375\\
5	0.3900146484375\\
7	0.3900146484375\\
9	0.3900146484375\\
11	0.3900146484375\\
13	0.3900146484375\\
15	0.3900146484375\\
17	0.3900146484375\\
19	0.3900146484375\\
21	0.3900146484375\\
23	0.3900146484375\\
25	0.3900146484375\\
27	0.3900146484375\\
29	0.3900146484375\\
};
\addlegendentry{ZF order 5}

\addplot [color=black, line width=1.0pt, draw=none, mark=asterisk, mark options={solid, black}]
  table[row sep=crcr]{%
1	0.4\\
3	0.4\\
5	0.4\\
7	0.4\\
9	0.4\\
11	0.4\\
13	0.4\\
15	0.4\\
17	0.4\\
19	0.4\\
21	0.4\\
23	0.4\\
25	0.4\\
27	0.4\\
29	0.4\\
};
\addlegendentry{Example fields}

\end{axis}
\end{tikzpicture}%
    	\caption{Robustness against different fields $f\in \mathcal{S}(1,L)$ for LPV system \eqref{eq:lpv_example}} 
    	\centering
    	\label{fig:LPV_generic_vehicle}
    \end{figure}

We now consider an LPV system $G(\rho)$ described by,
\begin{equation}\label{eq:lpv_example}
\begin{split}
\dot{x}&=v\\
\dot{v}&=-\rho(t)v -u,
\end{split}
\end{equation}
where, $\rho(t) \in \mathcal{P}=[0.8,1.2] \quad \forall t \in [0,\infty)$.
The scheduling parameter $\rho$ can be seen as a time-varying or adaptive damping co-efficient and can be either fixed and unknown or time-varying. 
It can be verified that \eqref{eq:perf_LMI_lpv} for this example is affine in $\Bar{\rho}$ and hence satisfaction of the inequality for $\Bar{\rho}=0.8$ and $\Bar{\rho}=1.2$ implies the satisfaction for any $\Bar{\rho}\in [0.8,1.2]$ \cite{.2000scherer}.
This reduces the condition \eqref{eq:perf_LMI_lpv} to a finite dimensional feasibility problem that is implemented to produce the results discussed next.
Fig. \ref{fig:LPV_generic_vehicle} shows the convergence rate estimates provied by different multipliers for fields $f \in \mathcal{S}(1,L)$ with increasing $L$.
As in the previous examples, example fields (quadratic) are chosen to get an upper bound on the convergence rate.
The reduction in conservatism with increasing order of multiplier can be clearly seen.
For this example, fifth order ZF multipliers show the least conservatism (possibly no conservatism).
This analysis essentially guarantees, that for this chosen example, poorly conditioned fields do not affect the convergence rate.
For constant parameter trajectories, i.e.,  $\rho(t)=\Bar{\rho} \in [0.8,1.2] \quad \forall t$, and for quadratic fields (linear gradients), a root-locus argument can be used to show that $\alpha=0.4$ for any $L>m$.
With the current performance analysis, we can see that this holds even for any non-constant trajectories $\rho$ restricted to the allowable parameter range and for any strongly convex field $f \in \mathcal{S}(m,L)$.
This example also illustrates the benefit of non-causal multipliers over causal ones and the reduction in conservatism with increasing order of the ZF multiplier.
\subsection{Quadrotor with two modes}
\begin{figure}[t]
%
%
\definecolor{mycolor1}{rgb}{0.00000,1.00000,1.00000}%
\definecolor{mycolor2}{rgb}{0.00000,0.44700,0.74100}%
\definecolor{mycolor3}{rgb}{0.85000,0.32500,0.09800}%
\definecolor{mycolor4}{rgb}{0.92900,0.69400,0.12500}%
\definecolor{mycolor5}{rgb}{0.49400,0.18400,0.55600}%
\begin{tikzpicture}

\begin{axis}[%
width=2.521in,
height=2.366in,
at={(0.758in,0.481in)},
scale only axis,
xmin=1,
xmax=7,
xlabel style={font=\color{white!15!black}},
xlabel={L},
ymin=0,
ymax=0.5,
ylabel style={font=\color{white!15!black}},
ylabel={$\alpha$},
axis background/.style={fill=white},
title style={font=\bfseries},
title={},
axis x line*=bottom,
axis y line*=left,
legend style={legend cell align=left, align=left, draw=white!15!black}
]
\addplot [color=green, dashed, line width=1.0pt]
  table[row sep=crcr]{%
1	0.3436279296875\\
1.5	0.2813720703125\\
2	0.234375\\
2.5	0.191650390625\\
3	0.150146484375\\
3.5	0.1092529296875\\
4	0.07080078125\\
4.5	0.0341796875\\
5	-1\\
5.5	-1\\
6	-1\\
6.5	-1\\
7	-1\\
7.5	-1\\
8	-1\\
8.5	-1\\
9	-1\\
9.5	-1\\
10	-1\\
};
\addlegendentry{$\text{ZF causal (P}_\text{1}\text{=0)}$}

\addplot [color=blue, dashdotted, line width=1.0pt]
  table[row sep=crcr]{%
1	0.3436279296875\\
1.5	0.2191162109375\\
2	0.1300048828125\\
2.5	0.0628662109375\\
3	0.00732421875\\
3.5	-1\\
4	-1\\
4.5	-1\\
5	-1\\
5.5	-1\\
6	-1\\
6.5	-1\\
7	-1\\
7.5	-1\\
8	-1\\
8.5	-1\\
9	-1\\
9.5	-1\\
10	-1\\
};
\addlegendentry{$\text{ZF anti-causal (P}_\text{3}\text{=0)}$}

\addplot [color=mycolor1, line width=1.0pt]
  table[row sep=crcr]{%
1	0.3436279296875\\
1.5	0.3070068359375\\
2	0.2593994140625\\
2.5	0.2117919921875\\
3	0.1654052734375\\
3.5	0.1214599609375\\
4	0.0799560546875\\
4.5	0.04150390625\\
5	0.0042724609375\\
5.5	-1\\
6	-1\\
6.5	-1\\
7	-1\\
7.5	-1\\
8	-1\\
8.5	-1\\
9	-1\\
9.5	-1\\
10	-1\\
};
\addlegendentry{ZF order 1}

\addplot [color=mycolor2, dotted, line width=1.0pt]
  table[row sep=crcr]{%
1	0.3436279296875\\
1.5	0.3070068359375\\
2	0.26123046875\\
2.5	0.2178955078125\\
3	0.177001953125\\
3.5	0.13916015625\\
4	0.1031494140625\\
4.5	0.069580078125\\
5	0.037841796875\\
5.5	0.0079345703125\\
6	-1\\
6.5	-1\\
7	-1\\
7.5	-1\\
8	-1\\
8.5	-1\\
9	-1\\
9.5	-1\\
10	-1\\
};
\addlegendentry{ZF order 2}




\addplot [color=black, line width=1.0pt, draw=none, mark=asterisk, mark options={solid, black}]
  table[row sep=crcr]{%
1	0.359877781991861\\
1.5	0.30771931207812\\
2	0.262397425883275\\
2.5	0.222141826976739\\
3	0.185818152359666\\
3.5	0.152649297281869\\
4	0.12207640065513\\
4.5	0.0936823881608809\\
5	0.0671467617957762\\
5.5	0.0422173233107017\\
6	0.0186916779750199\\
6.5	-0.00359531647237843\\
7	-0.0247803316318347\\
7.5	-0.0449778747672073\\
8	-0.0642849197875841\\
8.5	-0.0827843746763419\\
9	-0.100547722380696\\
9.5	-0.117637062795108\\
10	-0.13410671307487\\
};
\addlegendentry{Example fields}

\end{axis}
\end{tikzpicture}%
	\caption{Robustness against different fields $f\in \mathcal{S}(1,L)$ for a quadrotor with uncertain or time-varying mass $m\in \{0.2,2\}$} 
	\centering
	\label{fig:uncertain_quadrotor}
\end{figure}
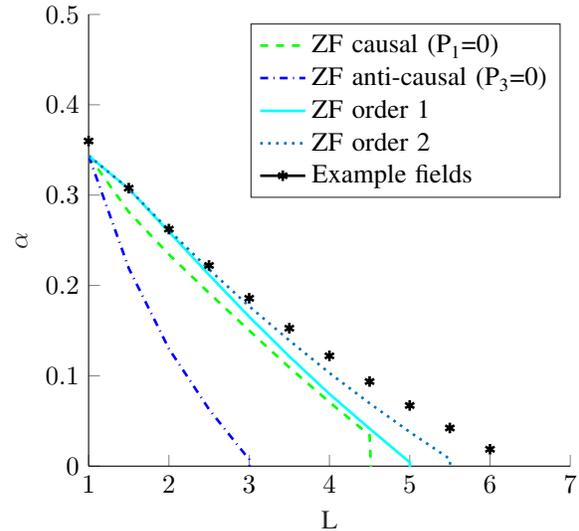

We now consider a scenario with a quadrotor, as in section \ref{eg:quadrotor}, but with two operating modes.
One operating mode corresponds to the quadrotor carrying some load and the other mode corresponds to no-load.
We model this by considering two masses $m \in \{0.2,2\}$ with LQR controllers designed as in section \ref{eg:quadrotor} for each mode separately.
We consider an arbitrary switching between the two modes and can be modeled as an LPV (or switching) system with $\mathcal{P}=\{1,2\}$ and $\rho(t) \in \mathcal{P} \quad \forall t$.
Fig. \ref{fig:uncertain_quadrotor} shows the convergence rate estimates provided by different multipliers for fields $f \in \mathcal{S}(1,L)$.
We observe that in comparison to the LTI case (Fig. \ref{fig:quadrotor_robustness} from section \ref{eg:quadrotor}), the performance is slightly reduced due to the possibility of arbitrary switching between modes.
Furthermore, the estimates with first order ZF multipliers are not tight anymore and we obtain better results with second order ZF multipliers.
No improvement in the estimates was observed upto 5th order ZF multipliers.
\section{Conclusions and Future Work}\label{sec:conclusions}
We have presented an approach to analyze the robust performance of source-seeking dynamics using the framework of $\alpha$-IQCs for LTI and LPV systems. By parameterizing the ZF multipliers, we have first derived a quasi-convex feasibility problem which can be used to obtain estimates on convergence rates.
We demonstrated the results on an LTI quadrotor model, an LTI non-minimum phase example where the estimates turned out to be tight with first order ZF multipliers and two LPV examples where the conservatism could be reduced by increasing the order of the ZF multiplier.
The reduction in conservatism by searching over general non-causal multipliers was also evident in some of the examples considered.
A number of interesting directions are open for further study.
Extension of these analysis results to controller synthesis is a valuable direction to be pursued.
Since convexity (not differentiability) of the underlying field is the key requirement in these results, source-seeking algorithms using sub-gradient forcing terms could be readily analyzed with minor changes in the theory.
Drawing motivation from \cite{hu2017unified}, an extension of these results to the stochastic setting is an interesting direction especially when considering imperfect communications.
Extending the full-block ZF multipliers from the standard literature on IQCs \cite{Fetzer.2017} to the setting of $\alpha-$IQCs is ongoing work.
On a different note, performance analysis of formation control algorithms including a single leader with gradient information has shown promising results and is currently investigated.




\section*{APPENDIX}
\subsection*{A.State-space realization of $\Psi$}
\begin{equation}\label{eq:ss_Psi}
\begin{split}
 &\left[
\begin{array}{c|c}
A_{\Psi} & B_{\Psi} \\
\hline
C_{\Psi} & D_{\Psi}
\end{array}
\right]\\
&=\\
&\left[
\begin{array}{cc|cc}
A_{\nu}^{\alpha} \otimes I_d        &\mathbf{0}                         &-B_{\nu}\otimes mI_d       &B_{\nu}\otimes I_d \\
\mathbf{0}                          &A_{\nu}^{\alpha} \otimes I_d       &B_{\nu}\otimes LI_d        &-B_{\nu}\otimes I_d \\
\hline
\mathbf{0}                          &\mathbf{0}                         &-mI_d                      &I_d \\
I_{\nu}\otimes I_d                  &\mathbf{0}                         &\mathbf{0}                 &\mathbf{0}\\
\mathbf{0}                          &\mathbf{0}                         &LI_d                       &-I_d\\
\mathbf{0}                          &I_{\nu}\otimes I_d                 &\mathbf{0}                 &\mathbf{0}
\end{array}
\right]   
\end{split}
\end{equation}

\subsection*{B.Proofs}
\begin{proof}[Theorem \ref{theom:lemma_pq_ZF}]
We proceed along the lines of \cite{Lessard.2016}.
Let $g: \mathbb{R}^d \xrightarrow[]{} \mathbb{R}$ be defined for any $\Tilde{y} \in \mathbb{R}^d$ by
\begin{equation}
    g(\Tilde{y})=f(\Tilde{y}+y_*)-f(y_*)-\frac{m}{2}||\Tilde{y}||^2,
\end{equation}
where, $f$ satisfies Assumption 1.
Assumption 1 can be used to show that $g \in \mathcal{S}(0,L-m)$, $g(\mathbf{0})=0$ and $\nabla g (\mathbf{0})=\mathbf{0}$. 
It can be shown \cite{Lessard.2016} that for all $\Tilde{y},\Tilde{y}_1,\Tilde{y}_2 \in \mathbb{R}^d$ , 
\begin{equation}\label{eq:signal_defn_r}
    (L-m)g(\Tilde{y}) -\frac{1}{2}||\nabla g(\Tilde{y})||^2 \geq 0
\end{equation}
\begin{equation}\label{eq:prop_g}
   (L-m)\nabla g(\Tilde{y})^T\Tilde{y} \geq (L-m)g(\Tilde{y}) +\frac{1}{2}||\nabla g(\Tilde{y})||^2
\end{equation}
\begin{equation}\label{eq:prop_g2}
\begin{split}
\nabla g(\Tilde{y}_1)^T&(\Tilde{y}_1-\Tilde{y}_2) \\
& \geq  g(\Tilde{y}_1)-g(\Tilde{y}_2)+ \frac{||\nabla g (\Tilde{y}_1) - \nabla g (\Tilde{y}_2)||^2}{2(L-m)}.
\end{split}
\end{equation}
Using \eqref{eq:signal_defn_r}, we define a non-negative function $r: \mathbb{R}^d \xrightarrow[]{} \mathbb{R}$ as,
\begin{equation}
    r(\Tilde{y})=(L-m)g(\Tilde{y}) -\frac{1}{2}||\nabla g(\Tilde{y})||^2.
\end{equation}
Using definitions \eqref{eq:signal_defns}, we verify that for $\Tilde{u},\Tilde{y}$ satisfying \eqref{eq:delta},
\begin{equation}
\begin{split}
    p(t) &= \nabla g(\Tilde{y}(t)),\\
    q(t) &= (L-m)\Tilde{y}(t)-\nabla g(\Tilde{y}(t)).
\end{split}
\end{equation}
Consider for $t_1 \in [0,\infty)$,
\begin{equation}\label{eq:pq_lower_bound1}
\begin{split}
p(t_1)^Tq(t_1)&=\nabla g (\Tilde{y}(t_1))^T ((L-m)\Tilde{y}(t_1)-\nabla g(\Tilde{y}(t_1)))\\
&=(L-m)\nabla g(\Tilde{y}(t_1))^T\Tilde{y}(t_1)-||\nabla g(\Tilde{y}(t_1))||^2\\
&\geq^{(\textnormal{using \ref{eq:prop_g}}) }(L-m)g(\Tilde{y}(t_1)) -\frac{1}{2}||\nabla g(\Tilde{y}(t_1))||^2\\
&=r(\Tilde{y}(t_1))\geq 0.
\end{split}
\end{equation}
We now consider signal extensions $\Tilde{u}_T,\Tilde{y}_T,p_T$ and $q_T$ as defined in \eqref{eq:extension_signal} and since the map \eqref{eq:delta} is static (and therefore causal), we have $\forall t \in \mathbb{R}$,
\begin{equation}\label{eq:signal_defns_estensions}
    \begin{split}
        p_T(t)&=\Tilde{u}_T(t)-m \Tilde{y}_T(t)\\
        q_T(t)&=L \Tilde{y}_T(t) - \Tilde{u}_T(t).
    \end{split}
\end{equation}

So, for $t_1,t_2 \in \mathbb{R}$,
\begin{equation}\label{eq:eq:pq_lower_bound2}
\begin{split}
p_T(t_1)^T[q_T(t_1)-&q_T(t_2)]\\
&=\\
\nabla g (\Tilde{y}_T(t_1))^T ( (L-m)&\Tilde{y}_T(t_1)-\nabla g(\Tilde{y}_T(t_1)) \\
 - (L-m)\Tilde{y}_T(t_2) &+\nabla g(\Tilde{y}_T(t_2)))\\
&=\\
(L-m)\nabla g (\Tilde{y}_T(t_1))^T(\Tilde{y}_T(t_1)-&\Tilde{y}_T(t_2)) - ||\nabla g (\Tilde{y}_T(t_1))||^2\\
+\nabla g (\Tilde{y}_T(t_1))^T& \nabla g (\Tilde{y}_T(t_2))\\
&\geq^{(\textnormal{using \ref{eq:prop_g2}})}\\
(L-m)(g(\Tilde{y}_T(t_1))-&g(\Tilde{y}_T(t_2)))+\\
\frac{1}{2}||\nabla g (\Tilde{y}_T(t_1))-
\nabla g (\Tilde{y}_T(t_2))||^2 &-||\nabla g (\Tilde{y}_T(t_1))||^2\\
+\nabla g (\Tilde{y}_T(t_1))^T & \nabla g (\Tilde{y}_T(t_2))\\
&=\\
((L-m)g(\Tilde{y}_T(t_1))-&\frac{1}{2}||\nabla g (\Tilde{y}_T(t_1))||^2)\\
-((L-m)g(\Tilde{y}_T(t_2))-&\frac{1}{2}||\nabla g (\Tilde{y}_T(t_2))||^2)\\
&=\\
r(\Tilde{y}_T(t_1))&-r(\Tilde{y}_T(t_2)).
\end{split}
\end{equation}
For any $\beta \in [0,1]$, multiplying \eqref{eq:pq_lower_bound1} by $(1-\beta)$, multiplying \eqref{eq:eq:pq_lower_bound2} by $\beta$ and adding, we get,
\begin{equation}\label{eq:pq_prop_lower_bound3}
    p_T(t_1)^T[q_T(t_1)-\beta q_T(t_2)] \geq r(\Tilde{y}_T(t_1))-\beta r(\Tilde{y}_T(t_2)).
\end{equation}
For any $\tau \in \mathbb{R}$, let $\beta(\tau)=\textnormal{min}\{1,e^{-2\alpha \tau}\}$. 
Noting that $\beta(\tau) \in [0,1]$ and $\beta(\tau)\leq e^{-2\alpha \tau} \quad \forall \tau \in \mathbb{R}$, we use \eqref{eq:pq_prop_lower_bound3} and non-negativity of $r$ to obtain,
\begin{equation}
\begin{split}
    p_T(t)^T[q_T(t)-&\beta(\tau)q_T(t-\tau)]\\
    &\geq \\
    r(\Tilde{y}_T(t))-&\beta(\tau) r(\Tilde{y}_T(t-\tau))\\
    &\geq \\
    r(\Tilde{y}_T(t))-&e^{-2\alpha \tau} r(\Tilde{y}_T(t-\tau)).\\
\end{split}
\end{equation}
Multiplying both sides by $e^{2\alpha t}$ and integrating from $0$ to $T$, we get,
\begin{equation}
\begin{split}
\int_0^T e^{2\alpha t}p_T(t)^T[q_T(t)&-\beta(\tau)q_T(t-\tau)] dt\\
&\geq\\
\int_0^T e^{2\alpha t} r(\Tilde{y}_T(t))-&e^{2\alpha(t-\tau)} r(\Tilde{y}_T(t-\tau)) dt\\
&=\\
\int_0^T e^{2\alpha t} r(\Tilde{y}_T(t)) dt -
& \int_{-\tau}^{T-\tau} e^{2\alpha s} r(\Tilde{y}_T(s)) ds.
\end{split}
\end{equation}
We finally use that the signal extensions are equal to the non-extended signals on $[0,T]$ to obtain,
\begin{equation}\label{eq:lower_bound_LHS_with_r}
\begin{split}
\int_0^T e^{2\alpha t}p(t)^T[q(t)&-\beta(\tau)q_T(t-\tau)] dt\\
& \geq \\
\int_0^T e^{2\alpha t} r(\Tilde{y}_T(t)) dt -
& \int_{-\tau}^{T-\tau} e^{2\alpha s} r(\Tilde{y}_T(s)) ds.
\end{split}    
\end{equation}
The integrand in both terms is identical and non-negative on $[0,T]$, and zero outside $[0,T]$. 
Hence if $\tau \geq 0$,
\begin{equation}
\int_{-\tau}^{T-\tau} e^{2\alpha s} r(\Tilde{y}_T(s)) ds =\int_{0}^{T-\tau} e^{2\alpha s} r(\Tilde{y}_T(s)) ds,
\end{equation}
and
\begin{equation}
\begin{split}
\int_0^T e^{2\alpha t} r(\Tilde{y}_T(t)) dt &-
\int_{-\tau}^{T-\tau} e^{2\alpha s} r(\Tilde{y}_T(s)) ds\\
&=\\
\int_{T-\tau}^T e^{2\alpha t} &r(\Tilde{y}_T(t)) dt \geq 0.
\end{split}
\end{equation}
For $\tau < 0$,
\begin{equation}
\int_{-\tau}^{T-\tau} e^{2\alpha s} r(\Tilde{y}_T(s)) ds =\int_{-\tau}^{T} e^{2\alpha s} r(\Tilde{y}_T(s)) ds,
\end{equation}
and
\begin{equation}
\begin{split}
\int_0^T e^{2\alpha t} r(\Tilde{y}_T(t)) dt &-
\int_{-\tau}^{T-\tau} e^{2\alpha s} r(\Tilde{y}_T(s)) ds\\
&=\\
\int_{0}^{-\tau} e^{2\alpha t} &r(\Tilde{y}_T(t)) dt \geq 0.
\end{split}
\end{equation}
This shows that the right hand side of \eqref{eq:lower_bound_LHS_with_r} is non-negative for all $\tau \in \mathbb{R}$ which completes the proof.
\end{proof}
\begin{proof}[Theorem \ref{theom:theorem_pqw_filtered_ZF}]
With the signal definitions \eqref{eq:signal_defns_w12}, an appropriate change of integration variable, signal extension as defined in \eqref{eq:extension_signal} and Fubini's theorem, we obtain,
\begin{align*}
    &\int_0^T e^{2\alpha t}p(t)^T w_1(t)dt \\
    &=\int_0^T e^{2\alpha t}p(t)^T \left( \int_{0}^t e^{-2\alpha (t-\tau)}h(t-\tau) q(\tau) d\tau \right)dt\\
    &=\int_0^T \int_{s=0}^t e^{2\alpha (t-s)}p(t)^T  h(s) q(t-s) ds dt \\
    &=\int_0^T \int_{s=0}^{\infty} e^{2\alpha (t-s)}p(t)^T  h(s) q_T(t-s) ds dt \\
    &=\int_{s=0}^{\infty} h(s) \left(\int_0^T  e^{2\alpha (t-s)}p(t)^T  q_T(t-s)dt \right)  ds. \numberthis  \label{eq:causal_term} 
\end{align*}
Similarly,
\begin{align*}
    &\int_0^T e^{2\alpha t}q(t)^T w_2(t)dt \\
    &=\int_0^T e^{2\alpha t}q(t)^T  \left(\int_0^t e^{-2\alpha (t-\tau)}h(-(t-\tau)) p(\tau) d\tau \right) dt\\
     &=\int_{t=0}^T \int_{\tau=0}^t e^{2\alpha \tau}q(t)^T h(\tau-t) p(\tau) d\tau dt\\
    &=\int_{\tau=0}^T \int_{t=\tau}^T  e^{2\alpha \tau}q(t)^T h(\tau-t) p(\tau) dt d\tau \\
    &=\int_{\tau=0}^T \int_{s=\tau-T}^0  e^{2\alpha \tau}q(\tau-s)^T h(s) p(\tau) ds d\tau \\
    &=\int_{\tau=0}^T \int_{s=-\infty}^0  e^{2\alpha \tau}q_T(\tau-s)^T h(s) p(\tau) ds d\tau \\
    &=\int_{s=-\infty}^0 h(s)  \left(\int_{\tau=0}^T   e^{2\alpha \tau}p(\tau)^Tq_T(\tau-s)  d\tau \right) ds \numberthis  \label{eq:anti-causal_term} 
\end{align*}

Putting \eqref{eq:causal_term} and \eqref{eq:anti-causal_term} together, we get,
\begin{equation}
\begin{split}
    &\int_0^T e^{2\alpha t}(p(t)^T w_1(t) + q(t)^T w_2(t)) dt \\
    &=\int_{s=0}^{\infty} h(s) \left( \int_0^T  e^{2\alpha (t-s)}p(t)^T  q_T(t-s)dt \right)ds \\
    &\quad \quad + \int_{s=-\infty}^0 h(s)  \left(\int_{\tau=0}^T   e^{2\alpha \tau}p(\tau)^Tq_T(\tau-s)  d\tau \right) ds\\
    &=\int_{s=-\infty}^{\infty} h(s)  \left(  \int_0^T  e^{2\alpha t} \textnormal{min}\{1,e^{-2\alpha s}\} p(t)^T   q_T(t-s)dt \right) ds \\
    & \leq^ {\textnormal{using $h(t)>0$ and Lemma \ref{theom:lemma_pq_ZF}}}\\
    &\int_{s=-\infty}^{\infty} h(s)  \left(  \int_0^T  e^{2\alpha t} p(t)^T   q(t)dt \right) ds \\
    & \leq ^ {\textnormal{using \eqref{eq:zf_impulse_resp_cond}}}\\
    &\int_0^T H e^{2\alpha t} p(t)^T   q(t)dt \\
\end{split}
\end{equation}
\end{proof}
\begin{proof}[Theorem \ref{theom:theorem_LMI_ZF}]
Let
\begin{equation}\label{eq:h_defn}
    \begin{split}
        h(t)&=
        \begin{cases}
        P_1 Q_{\nu}(-t) & \text{if}\ t<0, \\
        P_3 Q_{\nu}(t) & \text{if}\ t\geq 0.
        \end{cases}
    \end{split}
\end{equation}
It has been shown in \cite{J.Veenman.2013} that if $H, P_1, P_3$ satisfy \eqref{eq:L1_norm_constraint} and \eqref{eq:positivity}, $h$ defined in \eqref{eq:h_defn} satisfies \eqref{eq:zf_impulse_resp_cond}.
From the state-space realization \eqref{eq:ss_Psi} of $\Psi$ and the signal definitions \eqref{eq:signal_defns}, we get, 
\begin{equation}
\begin{split}
\Tilde{z}(t)&=
\begin{bmatrix}
\Tilde{u}(t)-m \Tilde{y}(t)\\
\int_0^t e^{A_{\nu}^{\alpha}(t-\tau)}B_{\nu}(\Tilde{u}(\tau)-m \Tilde{y}(\tau))d\tau\\
L \Tilde{y}(t) - \Tilde{u}(t)\\
\int_0^t e^{A_{\nu}^{\alpha}(t-\tau)}B_{\nu}(L \Tilde{y}(\tau) - \Tilde{u}(\tau))d\tau\\
\end{bmatrix}\\
&=\begin{bmatrix}
p(t)\\
\int_0^t e^{-2\alpha (t-\tau)}Q_{\nu}(t-\tau)p(\tau)d\tau\\
q(t)\\
\int_0^t e^{-2\alpha (t-\tau)}Q_{\nu}(t-\tau)q(\tau)d\tau\\
\end{bmatrix}.
\end{split}
\end{equation}
Using the block structure of matrices $P \in \mathbb{P}$,
\begin{equation}
\begin{split}
 &\Tilde{z}^T(t) (P \otimes I_d) \Tilde{z}(t) \\
 &= 2H p(t)^T q(t) \\
 & \quad \quad - 2p(t)^T\int_0^t e^{-2\alpha (t-\tau)}P_3 Q_{\nu}(t-\tau)q(\tau)d\tau \\
&\quad \quad -2q(t)^T\int_0^t e^{-2\alpha (t-\tau)}P_1 Q_{\nu}(t-\tau)p(\tau)d\tau.  
\end{split}
\end{equation}
Using the signal definitions \eqref{eq:signal_defns_w12} and \eqref{eq:h_defn}, we get,
\begin{equation}
\begin{split}
 \Tilde{z}^T(t) (P \otimes I_d) \Tilde{z}(t) =&2H p(t)^T q(t)-2p(t)^T w_1(t)\\
 &\quad \quad -2q(t)^Tw_2(t).  
\end{split}
\end{equation}
Theorem \ref{theom:theorem_pqw_filtered_ZF} can now be applied to finish the proof.
\end{proof}
\begin{proof}[Theorem \ref{theom:theorem_perf_analysis_lpv}]
Since $\left[\begin{array}{c|c}
\mathcal{A}(\rho)     &  \mathcal{B}(\rho)\\
\hline
\mathcal{C}(\rho)     &  \mathcal{D}(\rho)
\end{array}\right]$ represents the serial concatenation $\Psi \begin{bmatrix}
G(\rho)\\
I
\end{bmatrix}$,
let $\xi=\begin{bmatrix}
\Tilde{\eta} \\
x_{\psi}
\end{bmatrix}$ be the concatenation of the state $\Tilde{\eta}$ of $G(\rho)$ and the filter state $x_{\Psi}$.
For any trajectory $\rho$, such that, $\rho(t) \in \mathcal{P} \quad \forall t \in [0,\infty)$, the dynamics of $\xi$, can be represented by
\begin{equation} \label{eq:sys_dyn_psi_GI_lpv}
    \begin{split}
        \Dot{\xi}&=\mathcal{A}(\rho(t))\xi + \mathcal{B}(\rho(t))\Tilde{u}, \quad\quad \xi(0)=[\Tilde{\eta}_0^T \quad \mathbf{0}]^T,\\
        \Tilde{z}&=\mathcal{C}(\rho(t))\xi + \mathcal{D}(\rho(t))\Tilde{u}.
    \end{split}
\end{equation}
Since \eqref{eq:sys_dyn_psi_GI_lpv} is the serial concatenation $\Psi \begin{bmatrix}
G(\rho)\\
I
\end{bmatrix}$, 
the output $\Tilde{z}$ of \eqref{eq:sys_dyn_psi_GI_lpv} is the one defined in \eqref{eq:signal_defn_z} with $\Tilde{y}$ defined as the output of $G(\rho)$ for input $\Tilde{u}$.
Furthermore, dynamics  \eqref{eq:sys_dyn_G_lpv} imply that $\Tilde{u},\Tilde{y}$ satisfy \eqref{eq:delta}. 
Hence, Theorem \ref{theom:theorem_LMI_ZF} implies
\begin{equation} \label{eq:zPz_pos}
    \int_0^T e^{2\alpha t}\Tilde{z}^T(t) (P \otimes I_d) \Tilde{z}(t) dt \geq 0 \quad \forall P \in \mathbb{P},\forall T \geq 0.
\end{equation}

Define a storage function $V(\xi)=\xi^T \mathcal{X} \xi$.
Using \eqref{eq:sys_dyn_psi_GI_lpv}, \eqref{eq:perf_LMI_lpv} and the assumption that $\rho(t) \in \mathcal{P} \quad \forall t \in [0,\infty)$, we get,
\begin{equation*}
\begin{split}
&\frac{d}{dt}(V(\xi(t)))+2\alpha V(\xi(t))=\\
    &\begin{bmatrix}
    \xi & \Tilde{u}
    \end{bmatrix}
    \begin{bmatrix}
    \mathcal{A}(\rho(t))^T\mathcal{X}+\mathcal{X}\mathcal{A}(\rho(t))+2\alpha \mathcal{X}  & \mathcal{X}\mathcal{B}(\rho(t)) \\
    \mathcal{B}(\rho(t))^T\mathcal{X}    & \mathbf{0}
    \end{bmatrix}
    \begin{bmatrix}
    \xi\\ \Tilde{u}
    \end{bmatrix}\\ 
    &\leq
    -\begin{bmatrix}
    \xi & \Tilde{u}
    \end{bmatrix}
    \begin{bmatrix}
\mathcal{C}(\rho(t))^T \\
\mathcal{D}(\rho(t))^T 
\end{bmatrix}
( P\otimes I_d)
\begin{bmatrix}
\mathcal{C}(\rho(t)) & \mathcal{D}(\rho(t)) 
\end{bmatrix}
\begin{bmatrix}
    \xi\\ \Tilde{u}
    \end{bmatrix}\\
    &=-\Tilde{z}^T(t) (P \otimes I_d) \Tilde{z}(t)   
\end{split}.
\end{equation*}
Rearranging, multiplying by $e^{2\alpha t}$ and integrating from $0$ to $T$, we obtain,
\begin{equation*}\label{eq:diff_diss_ineq}
\begin{split}
\frac{d}{dt}(e^{2\alpha t}V(\xi(t))) + e^{2\alpha t}\Tilde{z}^T(t) (P \otimes I_d) \Tilde{z}(t) &\leq 0\\
V(\xi(T)) + \int_0^T e^{2\alpha t}\Tilde{z}^T(t) (P \otimes I_d) \Tilde{z}(t) dt &\leq e^{-2\alpha t}V(\xi(0))
\end{split}
\end{equation*}
Using \eqref{eq:zPz_pos} and $\mathcal{X}>0$,
\begin{equation*}
    \begin{split}
        V(\xi(T)) &\leq e^{-2\alpha t}V(\xi(0))\\
        & \implies \\
        ||\xi|| &\leq \sqrt{\textnormal{cond}(\mathcal{X})} \cdot ||\xi(0)|| \cdot  e^{-\alpha t}\\
        & \implies \\
        ||\Tilde{\eta}|| &\leq \sqrt{\textnormal{cond}(\mathcal{X})} \cdot ||\xi(0)|| \cdot  e^{-\alpha t}\\
        & \implies \\
        ||\Tilde{y}|| &\leq ||C_G|| \cdot ||\Tilde{\eta}|| \\
        & \leq ||C_G|| \sqrt{\textnormal{cond}(\mathcal{X})} \cdot ||\xi(0)|| \cdot  e^{-\alpha t} \\
        &= \kappa e^{-\alpha t}
    \end{split}
\end{equation*}
\end{proof}

%
%


\bibliographystyle{IEEEtran}
%
\bibliography{root}

\begin{thebibliography}{10}
\providecommand{\url}[1]{#1}
\csname url@samestyle\endcsname
\providecommand{\newblock}{\relax}
\providecommand{\bibinfo}[2]{#2}
\providecommand{\BIBentrySTDinterwordspacing}{\spaceskip=0pt\relax}
\providecommand{\BIBentryALTinterwordstretchfactor}{4}
\providecommand{\BIBentryALTinterwordspacing}{\spaceskip=\fontdimen2\font plus
\BIBentryALTinterwordstretchfactor\fontdimen3\font minus
  \fontdimen4\font\relax}
\providecommand{\BIBforeignlanguage}[2]{{%
\expandafter\ifx\csname l@#1\endcsname\relax
\typeout{** WARNING: IEEEtran.bst: No hyphenation pattern has been}%
\typeout{** loaded for the language `#1'. Using the pattern for}%
\typeout{** the default language instead.}%
\else
\language=\csname l@#1\endcsname
\fi
#2}}
\providecommand{\BIBdecl}{\relax}
\BIBdecl

\bibitem{senga2007development}
H.~Senga, N.~Kato, A.~Ito, H.~Niou, M.~Yoshie, I.~Fujita, K.~Igarashi, and
  E.~Okuyama, ``Development of spilled oil tracking autonomous buoy system,''
  in \emph{OCEANS 2007}.\hskip 1em plus 0.5em minus 0.4em\relax IEEE, 2007, pp.
  1--10.

\bibitem{khong2014multi}
S.~Z. Khong, Y.~Tan, C.~Manzie, and D.~Ne{\v{s}}i{\'c}, ``Multi-agent source
  seeking via discrete-time extremum seeking control,'' \emph{Automatica},
  vol.~50, no.~9, pp. 2312--2320, 2014.

\bibitem{ogren2004cooperative}
P.~Ogren, E.~Fiorelli, and N.~E. Leonard, ``Cooperative control of mobile
  sensor networks: Adaptive gradient climbing in a distributed environment,''
  \emph{IEEE Transactions on Automatic control}, vol.~49, no.~8, pp.
  1292--1302, 2004.

\bibitem{Datar.51220205152020}
A.~Datar, P.~Paulsen, and H.~Werner, ``Flocking towards the source: Indoor
  experiments with quadrotors,'' in \emph{2020 European Control Conference
  (ECC)}.\hskip 1em plus 0.5em minus 0.4em\relax IEEE, 5/12/2020 - 5/15/2020,
  pp. 1638--1643.

\bibitem{megretski1997system}
A.~Megretski and A.~Rantzer, ``System analysis via integral quadratic
  constraints,'' \emph{IEEE Transactions on Automatic Control}, vol.~42, no.~6,
  pp. 819--830, 1997.

\bibitem{J.Veenman.2013}
{J. Veenman} and {C. W. Scherer}, ``Stability analysis with integral quadratic
  constraints: A dissipativity based proof,'' in \emph{52nd IEEE Conference on
  Decision and Control}, 2013, pp. 3770--3775.

\bibitem{UlfJonsson.}
U.~Jönsson, ``Lecture notes on integral quadratic constraints,'' 2001.

\bibitem{Scherer.16052021}
\BIBentryALTinterwordspacing
C.~Scherer, ``Dissipativity and integral quadratic constraints, tailored
  computational robustness tests for complex interconnections.'' [Online].
  Available: \url{https://arxiv.org/pdf/2105.07401}
\BIBentrySTDinterwordspacing

\bibitem{Attallah.2020}
A.~Attallah, A.~Datar, and H.~Werner, ``Flocking of linear parameter varying
  agents: Source seeking application with underwater vehicles,''
  \emph{IFAC-PapersOnLine}, vol.~53, no.~2, pp. 7305--7311, 2020.

\bibitem{OlfatiSaber.2006}
R.~Olfati-Saber, ``Flocking for multi-agent dynamic systems: Algorithms and
  theory,'' \emph{IEEE Transactions on Automatic Control}, vol.~51, no.~3, pp.
  401--420, 2006.

\bibitem{Lessard.2016}
L.~Lessard, B.~Recht, and A.~Packard, ``Analysis and design of optimization
  algorithms via integral quadratic constraints,'' \emph{SIAM Journal on
  Optimization}, vol.~26, no.~1, pp. 57--95, 2016.

\bibitem{Hu.2016}
B.~Hu and P.~Seiler, ``Exponential decay rate conditions for uncertain linear
  systems using integral quadratic constraints,'' \emph{IEEE Transactions on
  Automatic Control}, vol.~61, no.~11, pp. 3631--3637, 2016.

\bibitem{Zhang.25022019}
\BIBentryALTinterwordspacing
J.~Zhang, P.~Seiler, and J.~Carrasco, ``Noncausal fir zames-falb multiplier
  search for exponential convergence rate.'' [Online]. Available:
  \url{https://arxiv.org/pdf/1902.09473}
\BIBentrySTDinterwordspacing

\bibitem{Freeman.62720186292018}
R.~A. Freeman, ``Noncausal zames-falb multipliers for tighter estimates of
  exponential convergence rates,'' in \emph{2018 Annual American Control
  Conference (ACC)}.\hskip 1em plus 0.5em minus 0.4em\relax IEEE, 6/27/2018 -
  6/29/2018, pp. 2984--2989.

\bibitem{Fazlyab.2018}
M.~Fazlyab, A.~Ribeiro, M.~Morari, and V.~M. Preciado, ``Analysis of
  optimization algorithms via integral quadratic constraints: Nonstrongly
  convex problems,'' \emph{SIAM Journal on Optimization}, vol.~28, no.~3, pp.
  2654--2689, 2018.

\bibitem{shamma1992linear}
J.~S. Shamma and J.~R. Cloutier, ``A linear parameter varying approach to gain
  scheduled missile autopilot design,'' in \emph{1992 American Control
  Conference}.\hskip 1em plus 0.5em minus 0.4em\relax IEEE, 1992, pp.
  1317--1321.

\bibitem{pfifer2015robustness}
H.~Pfifer and P.~Seiler, ``Robustness analysis of linear parameter varying
  systems using integral quadratic constraints,'' \emph{International Journal
  of Robust and Nonlinear Control}, vol.~25, no.~15, pp. 2843--2864, 2015.

\bibitem{michalowsky2016extremum}
S.~Michalowsky and C.~Ebenbauer, ``Extremum control of linear systems based on
  output feedback,'' in \emph{2016 IEEE 55th Conference on Decision and Control
  (CDC)}.\hskip 1em plus 0.5em minus 0.4em\relax IEEE, 2016, pp. 2963--2968.

\bibitem{nelson2018integral}
Z.~E. Nelson and E.~Mallada, ``An integral quadratic constraint framework for
  real-time steady-state optimization of linear time-invariant systems,'' in
  \emph{2018 Annual American Control Conference (ACC)}.\hskip 1em plus 0.5em
  minus 0.4em\relax IEEE, 2018, pp. 597--603.

\bibitem{.2000scherer}
C.~Scherer and S.~Weiland, ``Linear matrix inequalities in control,''
  \emph{Lecture Notes, Dutch Institute for Systems and Control, Delft, The
  Netherlands}, vol.~3, no.~2, 2000.

\bibitem{ACC2022code}
\BIBentryALTinterwordspacing
A.~Datar and H.~Werner, \emph{Robust Performance Analysis of Source-Seeking
  Dynamics with Integral Quadratic Constraints}.\hskip 1em plus 0.5em minus
  0.4em\relax Zenodo, 2021. [Online]. Available:
  \url{https://doi.org/10.5281/zenodo.5564776}
\BIBentrySTDinterwordspacing

\bibitem{hu2017unified}
B.~Hu, P.~Seiler, and A.~Rantzer, ``A unified analysis of stochastic
  optimization methods using jump system theory and quadratic constraints,'' in
  \emph{Conference on Learning Theory}.\hskip 1em plus 0.5em minus 0.4em\relax
  PMLR, 2017, pp. 1157--1189.

\bibitem{Fetzer.2017}
M.~Fetzer and C.~W. Scherer, ``Full--block multipliers for repeated,
  slope--restricted scalar nonlinearities,'' \emph{International Journal of
  Robust and Nonlinear Control}, vol.~27, no.~17, pp. 3376--3411, 2017.

\end{thebibliography}

\end{document}